\magnification 1200
\parindent=0pt
\font\bfeins=cmbx12
\input amstex
\input amssym

\documentstyle{amsppt}

\define\abs{\par\vskip 0.3cm\goodbreak\noindent}
\define\FF{\par\vfill\eject}
\define\newl{\par\noindent}
\define\nl{\bigskip\item{}}
\define\snl{\smallskip\item{}}
\define\inspr #1{\parindent=20pt\bigskip\bf\item{#1}}
\define\iinspr #1{\parindent=27pt\bigskip\bf\item{#1}}
\define\einspr{\parindent=0cm\bigskip}
\define\ot{\otimes}
\define\trr{\triangleright}
\define\trl{\triangleleft}

\advance\voffset by -2.5truecm

\topmatter
\title\vskip -2truecm{\hfill DAMTP-98-13}
      \abs\abs
      {\bfeins Actions of Multiplier Hopf Algebras}
\endtitle
\author Bernhard Drabant, Alfons Van Daele and  Yinhuo Zhang
\endauthor
\thanks
E-mail: {\tt b.drabant\@damtp.cam.ac.uk}
\newl
\hbox to 1.5cm{\ \hfill}{\tt alfons.vandaele\@wis.kuleuven.ac.be}
\newl
\hbox to 1.5cm{\ \hfill}{\tt zhang\@davinci.wis.kuleuven.ac.be}
\endthanks
\address
{\it Drabant:} DAMTP, University of Cambridge, Silver Street,
Cambridge CB3 9EW, UK
\newl
{\it Van Daele and Zhang:} Department of Mathematics, KU Leuven, 
Celestijnenlaan 200B, B-3001 Heverlee, Belgium
\endaddress
\date
April 1997
\enddate
\abstract
For an action $\alpha$ of a group $G$ on an algebra $R$ (over $\Bbb C$),
the crossed product $R\times_\alpha G$ is the vector space of $R$-valued
functions with finite support in $G$, together with the twisted
convolution product given by
$$(\xi \eta)(p) = \sum_{q \in G} \xi(q) \alpha_q (\eta (q^{-1}p))$$
where $p\in G$. This construction has been extended to the theory of
Hopf algebras. Given an action of a Hopf algebra $A$ on an algebra $R$,
it is possible to make the tensor product $R\ot A$ into an algebra by
using a twisted product, involving the action. In this case, the algebra
is called the smash product and denoted by $R\# A$. In the group case,
the action $\alpha$ of $G$ on $R$ yields an action of the group algebra
$\Bbb C G$ as a Hopf algebra on $R$ and the crossed $R\times_\alpha G$
coincides with the smash product $R\# \Bbb C G$.
\newl
In this paper we extend the theory of actions of Hopf algebras to actions
of multiplier Hopf algebras. We also construct the smash product and we
obtain results very similar as in the original situation for Hopf
algebras.
\newl
The main result in the paper is a duality theorem for such actions.
We consider dual pairs of multiplier Hopf algebras to formulate this
duality theorem. We prove a result in the case of an algebraic quantum
group and its dual. The more general case is only stated and will be proven
in a separate paper on coactions. These duality theorems for actions are
substantial generalizations of the corresponding
theorem for Hopf algebras. Also
the techniques that are used here to prove this result are slightly
different and simpler.
\endabstract
\endtopmatter

\document
\baselineskip=14pt
\leftheadtext\nofrills{B.~Drabant, A.~Van Daele and Y.~Zhang}
\rightheadtext\nofrills{Actions of Multiplier Hopf Algebras}
%
%
%
%
\head 1. Introduction\endhead
\abs
Let $R$ be an algebra over $\Bbb C$ with identity and let $G$ be a group.
Assume that $\alpha$ is an action of $G$ on $R$
by means of automorphisms of $R$.
So, for every $p \in G$ we have an automorphism $\alpha_p$ of $R$ and we also
have $\alpha_{pq} (x) = \alpha_p (\alpha_q(x))$ for all $x \in R$ and $p,q
\in G$.  It follows that
$\alpha_e(x) = x$ for all $x \in R$
when $e$ is the identity of $G$.  There is a standard procedure to associate a
new algebra $R \times_\alpha G$ to such a triple $(R,G,\alpha)$.  The elements
in this algebra are $R$-valued functions with finite support in $G$.  This set
is made into a complex vector space in the obvious way and the product of two
such functions $\xi, \eta$ is a {\it twisted convolution product}, given by
  $$(\xi \eta)(p) = \sum_{q \in G} \xi(q) \alpha_q (\eta (q^{-1}p))$$
whenever $p \in G$.  This new algebra $R \times_\alpha G$ is called the crossed
product or semi-direct product of $R$ with the action $\alpha$ of $G$.
In the trivial case where $R = \Bbb C$, and so where $G$ acts trivially on
$\Bbb C$, we get of course the group algebra $\Bbb C G$ of
$G$ over $\Bbb C$.  If $R$ is non-trivial, but if the action is still trivial,
it is clear that the crossed product is the tensor product of $R$ with the
group algebra $\Bbb CG$.  If fact, if the action of $G$ on $R$ is inner, we
still get an algebra, isomorphic with this tensor product algebra.  When $R$
is itself the group algebra $\Bbb CH$ of a group $H$ and when $\alpha$ is
coming from an action of $G$ on $H$, then the crossed product is
the group algebra of the semi-direct product of $G$ and $H$.  In the
general case, this crossed product construction provides examples of
non-trivial algebras.  Moreover, some nice properties can be proven. One
of them is the duality for such crossed products (see e.g.\ [C-M] and
[Q]). We will say more about
this further in the introduction.
\abs
The above construction has been generalized to actions of Hopf algebras.  This
is done as follows.  Again let $R$ be an algebra over $\Bbb C$ with identity.
Let $A$ be a Hopf algebra over $\Bbb C$.  By an action of $A$ on $R$ is
meant a bilinear map $a \otimes x \rightarrow ax$ from $A \otimes R$ to $R$
such that $R$ is a left $A$-module and such that
$$a(xy) = \sum (a_{(1)}x)(a_{(2)}y)$$
for all $x,y \in R$ and $a \in A$ (we use the Sweedler convention). Here,
it is also natural to require that $1x = x$ for all $x \in R$ where $1$
is the identity in $A$.  Again, a new algebra can be constructed.  In
this case,
the tensor product $R \otimes A$ is the underlying vector space over $\Bbb C$
while the product is again a {\it twisted product} given by
$$(x \otimes a)(y \otimes b) = \sum x(a_{(1)}y) \otimes a_{(2)} b$$
for all $a,b \in A$ and $x,y \in R$.  In this setting, the new algebra
is called
the smash product (rather than the crossed product) and it is
denoted by $R \#A$.
(We will discuss the terminology later
in this introduction).  Also here, this construction provides interesting
examples and nice results can be proven (see e.g.\ [M2]).
\abs
When $G$ is a group, acting on an algebra $R$ and when $A$ is the Hopf algebra
$\Bbb C G$, then it is quite straightforward to see that $A$ acts on $R$ by
means of the map
$$f \otimes x \rightarrow \sum f(p) \alpha_p(x)$$
where $x \in R$ and $f$ is a complex function with finite support in $G$
(considered
as an element in $\Bbb C G$).  It is also easy to check that the crossed
product of $R$ by the action $\alpha$ of $G$, as defined before, coincides
with the smash product of $R$ with $A$ is this case.  So, indeed, the smash
product for an action of a Hopf algebra extends the notion of the crossed
(or semi-direct)
product for an action of a group.
\abs
One of the main results in the theory of smash products is the duality.  Let
us formulate this in its simplest form.  We consider a finite-dimensional
Hopf algebra $A$ and its dual Hopf algebra $A^\prime$.  As before, let $A$
act on $R$ and consider the smash product $R \# A$.  The dual action of
$A^\prime$ on $R \# A$ is now defined by
$$\omega  (x \# a) = \sum \omega (a_{(2)})\, x \# a_{(1)}$$
whenever $x \in R$, $a \in A$ and $\omega \in A^\prime$.  It is
possible to look
at the bismash product $(R \# A) \# A^\prime$ and it turns out that
this algebra
is isomorphic with the algebra $M_n(R)$ of $n \times n$ matrices over $R$
where $n$ is the dimension of $A$.  It is remarkable that the bismash product
no longer depends on the underlying Hopf algebra, nor on the action.  A similar
result still holds in more general situations (see [M2]).
\abs
This duality is an extension of the duality for crossed products of
actions of finite abelian groups. In this case, the dual action is an
action of the dual group. The duality is related with the Pontryagin
duality between a group $G$ and its dual group $\hat G$. This theory for
group actions has been extended in various directions. One of the
fundamental results has been obtained in the von Neumann algebra case.
There, a locally compact abelian group $G$ acts on a von Neumann algebra
$M$ in a continuous way and the crossed product $M\times _\alpha G$ is a
von Neumann algebra, constructed in a way very similar to the case
described above. The Pontryagin dual is now again an abelian locally
compact group $\hat G$ and the dual action $\hat \alpha$ is a natural
action of $\hat G$ on the crossed product $M\times_\alpha G$. The double
crossed product $(M\times_\alpha G)\times_{\hat \alpha} \hat G$ turns
out to be canonically isomorphic with the von Neumann algebra tensor
product $M\overline\ot \Cal B(\Cal L^2(G))$ of the original von Neumann
algebra $M$ with the algebra of all bounded linear operators on the
Hilbert space $\Cal L^2(G)$ of all square integrable functions w.r.t.\
a left Haar measure on $G$.
\abs
This duality for actions of locally compact groups on von Neumann
algebras was in fact the first duality of this kind that was obtained
and the duality for actions of Hopf algebras was inspired by this
result. It was proven by Takesaki [Ts] (see also [N-T] and [VD1]).
\abs
But this result has also been generalized within the field of operator
algebras in various directions. First, there is the duality for actions
of locally compact groups on C$^*$-algebras [T], see aslo [P].
There are also different
versions for non-abelian groups. Different aspects of non-abelian
abstract harmonic analysis come in. One of the more sophisticated
results is the duality for actions of Kac algebras on von Neumann
algebras [E-S]. In the operator algebra context, it seems more natural to work
with coactions, rather than with actions. In fact recently, some work
has been obtained for coactions of Hopf C$^*$-algebras [Ng].
\abs
All this work
points in the direction of {\it locally compact quantum groups}.
Unfortunately, there still is no satisfactory definition of such a
locally compact quantum group.
Recently, there has been some progress in developing such a notion. We
have the work of Masuda, Nakagami and Woronowicz ([M-N] and [M-N-W]).
There is also the
theory, developed by one of the authors on multiplier Hopf algebras and
algebraic quantum groups [VD5, VD6]. The C$^*$-version can be found in
[K-VD] and [K]. The theory of algebraic quantum groups is not so general as the
theory of Masuda and his coauthors. On the other hand, it is more
accessable and still it contains many interesting cases, such as the
discrete and the compact quantum groups. It also behaves well with
respect to duality.
\abs
Because of this, and because this theory of algebraic quantum groups is
close in nature to the theory of Hopf algebras, it is quite natural to
look at the theory of actions of Hopf algebras and to try to extend it
to the multiplier Hopf algebras. This is in fact the {\it content of this
paper}.
\abs
When a Hopf algebra $A$ is infinite-dimensional, the space $A^\prime$ of all
linear functionals is still an algebra, but no longer a Hopf algebra.
The natural candidate for the comultiplicaton on $A^\prime$ will have values in
$(A \otimes A)^\prime$ but this space is strictly bigger than $A^\prime \otimes
A^\prime$.  In some cases, this problem can be overcome by taking the subspace
$A^0$ of elements in $A^\prime$ for which the image of the
comultiplication lies
in $A^\prime \otimes A^\prime$.  This is a Hopf algebra, but unfortunately, in
general, it may be very small.  When it is big enough, we are essentially
dealing with a (non-degenerate) dual pair of Hopf algebras.  The
setting of dual
pairs is a natural framework to formulate results that are easy to obtain for
the finite-dimensional case (because $A^\prime$ is again a Hopf algebra), and
cannot be generalized to arbitrary Hopf algebras.
\abs
A typical such result is the construction of the quantum double by
Drinfel'd [D] (see also [VD2]).
Also the duality for actions is a result of this nature (see [M2]).
\abs
Recently, the theory of Hopf algebras has been generalized to the so-called
multiplier Hopf algebras and also the dual pairs of such multiplier Hopf
algebras have been studied [D-VD].  The quantum double could be
constructed for such
dual pairs.  And it should be no surprise that also the theory of actions and
the duality can be generalized to multiplier Hopf algebras and dual pairs of
multiplier Hopf algebras.
\abs
In section 2, we first recall the notion of a multiplier Hopf algebra.  Very
roughly speaking, it is the natural extension of the notion of a Hopf algebra
where the underlying algebra is no longer assumed to have an identity.
It has been argued
in previous papers (see e.g.\ [VD3] and [VD6])
why such an extension is important.  Apart from
recalling the definition of a multiplier Hopf algebra, we also prove some new
results for general multiplier Hopf algebras (existence of local units),
important as such, but especially nice for dealing with modules further in the
paper.  In section 2, we also recall some special cases of multiplier Hopf
algebras (with the existence of integrals and cointegrals).
\abs
In section 3, we treat modules.  Multiplier Hopf algebra modules where
already considered in [D-VD]. Here we investigate them closer.
We have special attention for the features that
are new because our algebras no longer are assumed to have an identity.  The
natural condition $1x = x$ for a module over an algebra with identity
will be replaced by the condition $AR = R$ (where $A$ is an algebra and
$R$ a left $A$-module).  We will explain where this condition comes
from and why it is natural.  We call such a module unital.
\abs
In section 4, we consider actions of a multiplier Hopf
algebra $A$ on an algebra $R$ and define the notion of an $A$-module
algebra.  To begin with, $R$ is a unital left $A$-module.
Moreover, we require that also
$$a(xy) = m(\Delta (a)(x \otimes y))$$
whenever $x,y \in R$ and $a \in A$ (and where $m$ denotes multiplication as a
linear map from $R \otimes R$ to $R$).  Precisely, because the module $R$ is
assumed to be unital, we can think of elements $y$ as elements of the form
$bz$ where $b \in A$ and $z \in R$ and therefore the above right hand side
can be read as
$$m((\Delta (a)(1 \otimes b))(x \otimes z))$$
and for multiplier Hopf algebras, we do have that $\Delta (A) (1 \otimes A)
\subseteq A \otimes A$.
\abs
In section 4, we also consider the adjoint action and related notions such as
inner actions, cocycle equivalent actions,... One can define the fixed
point algebra within $R$ but it turns out to be more natural first to
extend the action of $A$ to the multiplier algebra $M(R)$ and then to
consider fixed points in $M(R)$. Also this is treated in the fourth
section.
\abs
In section 5, we construct the smash product for an action of a multiplier Hopf
algebra $A$ on an algebra $R$.  It will be denoted by $R \# A$.  We obtain some
natural properties.  We get a universal property characterizing $R \# A$.  We
obtain that $R \# A$ is the tensor product algebra $R \otimes A$ when
the action
is inner and we find that cocycle equivalent actions yield isomorphic smash
products.
\abs
The main section is section 7 where we obtain our duality results.
However, we first need to study some special cases of the previous
theory. In all the previous sections, we have used the case of a group
action as a
motivating example.  In section 6 however, we treat the less trivial examples.
The basic framework is that of a dual pair of multiplier Hopf algebras as
introduced in [D-VD].  Moreover, we relate the results of this paper with the
constructions needed to obtain the quantum double.  As a special case, we have
the pairing of an algebraic quantum group $A$ with its dual $\hat A$
as in  [VD5, VD6].  This
case is closest to the duality of finite-dimensional Hopf algebras.  But
it is already far more general.
\abs
The duality theorem in section 7 is formulated in the same framework of
dual pairs. However, it turns out that the conditions we need to prove
this theorem for arbitrary dual pairs makes the result in fact a theorem
about coactions. Since we are preparing a separate paper on coactions of
multiplier Hopf algebras ([VD-Z3]), we will not give the proof of the
general result here. For algebraic quantum groups, the result is nicer
and easier to proof. The conditions that we need in the more general
case turn out to be automatic for algebraic quantum groups. So, the
duality theorem that we obtain here for algebraic quantum groups is
really a special case of the one obtained in [VD-Z3] for more general
pairings.
We will explain all this in this last section.
\abs
This more restricted duality (in the case of algebraic quantum
groups) turns out to be very similar to the original duality results
proven in the theory of actions of locally compact groups on von Neumann
algebras and C$^*$-algebras.
\abs
Let us now fix some  terminology, notations and  standard
references.
\abs
The standard references for Hopf algebras are [A] and [S]. A good
introduction can also be found in [M2]. For multiplier Hopf algebras and
algebraic quantum groups, we refer to [VD3, VD5, VD6]. See also [VD4] and
[VD-Z1] for related material. The C$^*$-version of algebraic quantum groups has
been treated in [K-VD] and [K]. For dual pairs of Hopf algebras, we have
[VD2] and for dual pairs of multiplier Hopf algebras (more important for
this paper), we have [D-VD]. For coactions, we have [VD-Z3].
\abs
For actions of Hopf algebras and for smash
products, see [M2]. Some standard results about actions of locally
compact groups on von Neumann algebras and C$^*$-algebras can be found
in [P] and [VD1].
\abs
In general, we work with (associative)
algebras over $\Bbb C$. Our algebras need not
have an identity, but the product is always required to be
non-degenerate (as a bilinear map). The fact that our algebras need not
have an identity is very important, but it also causes some extra
problems. Throughout the paper, we will explain on several occasions
what kind of problems occur. We will have a more or less standard way to
overcome these. This will be explained in the different situations.
If we have an identity (like in the multiplier algebra), we will
always denote it by $1$.
\abs
The space of all linear functionals on a vector space $A$ (the dual)
will always be denoted by $A'$.
Tensor products are algebraic tensor products (as we are not working in
a topological context). We identify $A\ot \Bbb C$ and $\Bbb C\ot A$ with
$A$. We use $\iota$ for the identity map. When $\omega$ is a linear
functional on $A$, the slice maps $\iota\ot\omega$ and $\omega\ot\iota$
are maps from $A\ot A$ to $A$.
\abs
We use $\Delta$ for the comultiplication, both for Hopf algebras and for
multiplier Hopf algebras. In both cases, we also use $S$ for the antipode
and $\epsilon$ for the counit. Ocassionally, we will use $m$ for the
multiplication in an algebra $A$ as a linear map from $A\ot A$ to $A$.
Sometimes however, we use $m$ to denote a multiplier.
\abs
We will work with the Sweedler summation convention and denote
$\Delta(a)$ by $\sum a_{(1)}\ot a_{(2)}$, not only for usual Hopf
algebras, but also for multiplier Hopf algebras. This presents no
problem in the first case, but it is somewhat problematic for multiplier
Hopf algebras. The use of the Sweedler notation for multiplier Hopf
algebras was introduced in [D-VD]. In section 2 of this paper, we will
give some more comments about this topic. There, we will also say more
about the use of the Sweedler notation in this paper.
\abs
Not only when using the Sweedler notation, but also in the case of other
summations, we will mostly omit the summation variable. This should be
clear in these cases.
\abs
Let us finish this introduction with a  remark.
\abs
It should be clear from the
discussion above that the theory of actions and crossed (or smash)
products has developed in different fields of mathematics. We have
actions of Hopf algebras on the one hand, but also all sorts of actions
in the theory of operator algebras, involving different aspects of
abstract harmonic analysis. This can also be seen from a view at the
list of references.
The terminology is not always the same. In operator algebras, one is
used to work with actions, covariant systems, representations, covariant
representations, crossed products, ... Crossed products and tensor
products are usually completed for an appropriate topology. In algebra,
one works with modules, module algebras, smash products, ... And of
course, there is no completion and tensor products are just algebraic
tensor products.
\abs
Also the three authors of this paper have
different {\it roots}. These different backgounds have led to a fruitful
cooperation.
In operator algebras e.g., one is used to work with algebras
without identity. Considering those has led to the notion of multiplier
Hopf algebras; in fact, a very natural extension of the notion of Hopf
algebras if the assumption is dropped that 'all algebras have an
identity'. This has also led to more cases where integrals exist (e.g.\
the algebraic quantum groups [VD6] and in particular, the discrete
quantum groups ([E-R] and [VD4])).
\abs
It also had some influence, not only on the terminology
used in this paper, but also on the kind of results that are proven.
\abs
We have tried to make this paper accessable for
readers with these different backgrounds (algebra, analysis,
mathematical physics, ...). In fact, we do hope that our theory will
eventually lead to a comprehensive theory of actions of locally compact
quantum groups, including not only the actions of Hopf algebras and multiplier
Hopf algebras, but also the different types of actions of locally
compact groups and its generalizations.
\abs
\abs
\head 2. Multiplier Hopf Algebras\endhead
\abs
Multiplier Hopf algebras were introduced in [VD3] in 1994. Since the
introduction, certain new results have been
obtained. In [VD6] we have given a characterization of the regular
multiplier Hopf algebras in terms of the bijectivity of the antipode. In
[D-VD] we have introduced the use of the Sweedler notation for
multiplier Hopf algebras.
\snl
In this section we first recall some of the main definitions and results in
the theory of multiplier Hopf algebras.
But also in this paper, we will obtain some new results (see proposition
2.2 and 2.6 below).
\snl
First, recall the definition of a multiplier Hopf algebra (see [VD3] for
details).
\snl
Let $A$ be an algebra over $\Bbb C$, with or without identity, but with a
non-degenerate product.  The {\it multiplier algebra} $M(A)$ can be
characterized as the largest algebra with identity in which $A$ sits as a
dense (= essential) two-sided ideal.  Consider the tensor product
$A \otimes A$ of $A$ with itself.  This is again an algebra over $\Bbb C$
with a non-degenerate product.  The imbedding of $A \otimes A$ in the
multiplier algebra $M(A \otimes A)$ is factored through the natural
imbeddings of $A \otimes A$ in $M(A) \otimes M(A)$ and of
$M(A) \otimes M(A)$ in $M(A \otimes A)$.
\snl
A {\it comultiplication} (or a {\it coproduct}) on $A$ is a homomorphism
$\Delta : A \rightarrow M(A \otimes A)$ such that $\Delta (a) (1 \otimes b)$
and $(a \otimes 1)\Delta(b)$ are elements of $A \otimes A$ for all
$a,b \in A$.  It is assumed to be {\it coassociative} in the sense that
$$
(a \otimes 1 \otimes 1)(\Delta \otimes \iota)  (\Delta (b)(1 \otimes c))
= (\iota \otimes \Delta)((a \otimes 1) \Delta (b))(1 \otimes 1 \otimes c)
$$
for all $a,b,c \in A$ (where $\iota$ denotes the identity map).
\abs
Then recall the following definition (see [VD3]).

\inspr{2.1} Definition \rm
A pair $(A,\Delta)$ of an algebra $A$ over $\Bbb C$ with a
non-degenerate product and a comultiplication
$\Delta$ on $A$ is called a {\it
multiplier Hopf algebra} if the linear maps from $A \otimes A$ to itself,
defined by
$$\align
a \otimes b & \rightarrow \Delta (a)(1 \otimes b) \\
a \otimes b & \rightarrow (a \otimes 1) \Delta (b)
\endalign
$$
are bijective.
\einspr
The motivating example comes from a group.  If $G$ is any group, let $A$ be
the algebra of complex functions on $G$ with finite support and identify
$A \otimes A$ with the algebra of complex functions on $G \times G$ with
finite support and $M(A \otimes A)$ with the algebra of all complex functions.
Then $A$ is made into a multiplier Hopf algebra if we define $\Delta$ by
$(\Delta f)(p,q) = f(pq)$ whenever $f \in A$ and $p,q \in G$.
\snl
For any multiplier Hopf algebra, we have a {\it counit} and an {\it antipode}.
The counit is the unique linear map $\epsilon : A \rightarrow \Bbb C$
satisfying
$$\align
& (\epsilon \otimes \iota)(\Delta (a)(1 \otimes b)) = ab \\
& (\iota \otimes \epsilon)((a \otimes 1) \Delta (b)) = ab
\endalign
$$
for all $a,b \in A$. It is a homomorphism.  The antipode is the unique linear
map $S : A \rightarrow M(A)$ satisfying
$$\align
& m(S \otimes \iota)(\Delta (a)(1 \otimes b)) = \epsilon (a)b \\
& m(\iota \otimes S)((a \otimes 1)\Delta (b)) = \epsilon (b)a
\endalign
$$
where $m$ denotes multiplication defined as a linear map from $A \otimes A$ to
$A$ and extended to $M(A) \otimes A$ and $A \otimes M(A)$.  The antipode is a
anti-homomorphism.
\snl
In the group example, $\epsilon$ is given by $\epsilon (f) = f(e)$ where $e$
is the identity in $G$ and $S$ is given by $(S(f))(p) = f(p^{-1})$ whenever
$f \in A$ and $p \in G$.
\snl
If $A$ is an algebra with identity, then $M(A) = A$ and $M(A \otimes A) = A
\otimes A$ and clearly, when $\Delta$ is a comultiplication on $A$ making it
into a multiplier Hopf algebra, then it is actually a Hopf algebra.
Conversely, every Hopf algebra (over $\Bbb C$) is a multiplier Hopf algebra.  So, multiplier
Hopf algebras are the natural extensions of Hopf algebras to the case where
the underlying algebra has no identity.
\snl
If $\Delta$ is comultiplication on $A$ such that also $\Delta (a)(b\otimes 1)$
and $(1 \otimes a)\Delta(b)$ are in $A \otimes A$ for all $a,b \in A$, we can
consider the opposite comultiplication $\Delta^\prime$ obtained from $\Delta$
by composing it with the flip on $A \otimes A$ (extended to $M(A \otimes A))$.
A multiplier Hopf algebra $(A,\Delta)$ is called {\it regular} if also $(A,
\Delta^\prime)$ is a multiplier Hopf algebra (i.e.\ satisfies the conditions in
definition 2.1).  Regularity is automatic in certain special cases: when $A$ is
abelian, when $\Delta$ is coabelian and when $A$ is a $^\ast$-algebra
and $\Delta$
a $^\ast$-homomorphism.  The regular multiplier Hopf algebras are precisely
those for which the antipode maps $A$ into $A$ and is bijective (see [VD6]).  
In
the $^\ast$-algebra case, when $\Delta$ is a $^\ast$-homomorphism, the
antipode satisfies the relation $S(a^\ast) = S^{-1}(a)^\ast$ for all
$a \in A$.
\snl
In this paper we will only work with regular multiplier Hopf algebras.
Nevertheless, we have the feeling that many results are still valid for
multiplier Hopf algebras that are not regular. To obtain such results,
one has to be much more careful. For convenience,
we have chosen to concentrate first on the regular case.
\snl
So, in what follows let $(A,\Delta)$ be a regular multiplier Hopf algebra.

\inspr{2.2} Proposition \rm
Let $a_1,a_2,...,a_n$ be elements in $A$.  Then there exist
elements $e,f$ in $A$ such that $ea_i = a_i$ and $a_i f = a_i$ for all $i$.
\snl
\bf Proof : \rm  Take elements $a_1,a_2,...,a_n$ in $A$ and consider the
subspace $V$ of $A^n$ given by
$$V = \{ (aa_1,aa_2,...,aa_n) \mid  a  \in A \}.$$
The existence of such an element $e$ follows if
$(a_1,a_2,...,a_n) \in V$.  Suppose this is not the case.  Then, choose
linear functionals $\omega_1,\omega_2,...,\omega_n$ on $A$ such that
$$\sum \omega_i (aa_i) = 0$$
for all $a \in A$ but $\sum \omega_i (a_i) \ne 0$.
Take any $a \in A$ and write
$$\Delta (a)(1 \otimes a_i) = \sum_j p_{ij} \otimes q_{ij}$$
with $p_{ij}, q_{ij} \in A$.  For all $b \in A$, we get
$$\sum_{i,j} \omega_i  (bq_{ij}) p_{ij}
= \sum_i (\iota \otimes \omega_i)((1 \otimes b)\Delta(a)(1 \otimes a_i))
= 0.
$$
Now, take a basis $(e_k)$ in the (finite-dimensional) subspace $W$ of $A$
spanned by the elements $\{ Sp_{ij} \mid  i,j \}$ and consider a dual
basis, i.e.\
linear functionals $(f_k)$ such that $f_k(e_j) = 0$ when $k \ne j$ and $f_k(e_k)
= 1$.  Then, for all $x \in W$ we have $\sum_k f_k (x) e_k = x$.  Hence
$$\align
\sum_{i,j} \omega_i(S(p_{ij})q_{ij}) & = \sum_{i,j,k} \omega_i (e_k q_{ij})f_k
(S(p_{ij})) \\
& = \sum_k f_k (S(\sum_{i,j} \omega_i (e_k q_{ij})p_{ij})) = 0.
\endalign
$$
Now, for all $i$ we have $\sum_j S(p_{ij})q_{ij} = \epsilon(a)a_i$ by the
property of the antipode.  Hence $\sum \omega_i (a_i) \epsilon (a) = 0$.  So, if
we choose $a$ such that $\epsilon(a) \ne 0$ we would get $\sum \omega_i (a_i) =
0$ and this leads to a contradiction.
\snl
This proves the existence of $e$.  Similarly, or by applying the antipode, we
get the existence of $f$. \hfill $\blacksquare$
\einspr

As we can see in the proof, we do not need the fact that $S$ is bijective to
prove that $e$ exists, but we have used that $S(A) \subseteq A$.  So, it is not
clear if the previous result remains true for any multiplier Hopf algebra. (In
fact we even don't know of examples of multiplier Hopf algebras for which
$S(A)$ does not lie in $A$ !).
\snl
The property of $A$ that we proved in proposition 2.2 has some simple but nice
(and important) consequences.
\snl
First remark of course that the property is trivial for algebras with an
identity.
\snl
Perhaps the main consequence is the justification of the use of the Sweedler
notation, certainly for this paper.  The Sweedler notation for multiplier Hopf
algebras was first used in [D-VD]. Using the property above, we can justify the
Sweedler notation in this context even more.
\snl
We would like to use a formal expression for $\Delta (a)$ when $a \in A$.  The
problem is that $\Delta (a)$ is not in $A \otimes A$ in general.  We know
however that $\Delta(a)(1 \otimes b) \in A \otimes A$ for all $a,b \in A$.  We
have proposed in [D-VD] to write $\sum_{(a)} a_{(1)} \otimes a_{(2)} b$ for this.
Now, we know that there is an element $e \in A$ such that $b = eb$ and we can
think of $\sum_{(a)} a_{(1)} \otimes a_{(2)}$ to stand for $\Delta(a)(1
\otimes e)$.  Of course, this is still dependent on $b$.  But we know that
for several elements $b$, we can use the same $e$.
\snl
When using Sweedler's notation for regular multiplier Hopf algebras, one
always has to make sure that at most one factor $a_{(k)}$ is not {\it covered} by an
element of $A$.  Expressions like $\sum a_{(1)} \otimes ba_{(2)}$ and $\sum
a_{(1)} \otimes ba_{(2)} \otimes a_{(3)} c$ are allowed.  For this last
expression take e.g.\ elements $e,f$ in $A$ such that $bf = b$ and $ec = c$ and
use $\sum a_{(1)} \otimes fa_{(2)} \otimes a_{(3)}e$ for
$(1 \otimes f \otimes 1)(\Delta \otimes \iota)(\Delta (a)(1 \otimes e))$ which
is well-defined in $A \otimes A \otimes A$.
\snl
Nevertheless, the Sweedler notation should be used with some care and it is
good to have in mind what the expression stands for. One situation where
one has to be particularly careful is when using the Sweedler notation
in defining maps (because of the ambiguity).  If any doubt, it is
always possible to avoid it.  On the other hand, the Sweedler notation is
very convenient and makes formulas much more transparant.  It has been our
choice to make use of it in this paper, also if, at first sight, it is
somewhat dangerous. In any case, the reader can easily translate the formulas
to a form which does not contain the Sweedler notation.
\snl
Let us now turn our attention to a second consequence of proposition 2.2.
First recall that for any $a \in A$ and $\omega \in A^\prime$ we can define
elements $(\omega \otimes \iota)\Delta(a)$ and $(\iota \otimes \omega)\Delta
(a)$ in $M(A)$.  Now we can look at the subspace of $A^\prime$ of elements
$\omega$ such that $(\omega \otimes \iota)\Delta (a)$ and $(\iota \otimes
\omega) \Delta (a)$ are in $A$ for all $a \in A$.  We know that elements of
the form $\omega(\,\cdot\, b)$ and $\omega(b \,\cdot\,)$ where $\omega \in A^\prime$
and $b \in A$ have this property.  Now it follows from the proposition that
these functionals also form a subspace of $A^\prime$.  Indeed, take
e.g.\ $\omega, \psi \in A^\prime$ and $b,c \in A$.  Choose $e \in A$ so that
$eb = b$ and $ec = c$.  Define $\rho \in A^\prime$ by $\rho(a) = \omega (ab)
+ \psi(ac)$.  Then $\rho = \rho (\,\cdot\, e)$ so that also $\rho$ has this form.
\snl
It also follows from this proposition that elements of the form $\rho = \omega
(\,\cdot\, b)$ have unique extensions to $M(A)$ satisfying $\rho(x) = \omega (xb)$
for all $x \in M(A)$.  Indeed, if again $\omega, \psi \in A^\prime$ and
$b, c
\in A$ and if now $\omega (ab) =\psi (ac)$ for all $a \in A$, then again take
$e$ such that $eb = b$ and $ec = c$ and we find for all $x \in M(A)$ that
$\omega (xb) = \omega (x eb) = \psi (xec) = \psi (xc)$.
\snl
These two other consequences are not so important for this paper. They
were used e.g.\ in [VD-Z1].
\nl
Next, let us recall the notion of an integral.

\inspr{2.3} Definition \rm
Let $(A, \Delta)$ be a multiplier Hopf algebra.  A linear functional $\varphi$
is called a {\it left integral} if it is non-zero and if $(\iota \otimes
\varphi)\Delta (a) = \varphi(a)1$ for all $a \in A$.  A non-zero linear
functional $\psi$ is called a {\it right integral} if $(\psi \otimes
\iota)\Delta (a) = \psi (a)1$ for all $a \in A$.
\einspr

It has been shown in [VD6] that integrals are unique (up to a scalar of course)
for regular multiplier Hopf algebras.  They are also faithful (i.e.\ if $\varphi
(a\,\cdot\,) = 0$ then $a = 0$ and if $\varphi(\,\cdot\, a) = 0$, then $a = 0$).
They  satisfy the so-called weak K.M.S. property: there is an automorphism
$\sigma$ of $A$ such that $\varphi(ab) = \varphi (b \sigma (a))$ for all $a,b
\in A$.
It is also clear that the antipode converts a left integral to a right one and
the other way around.
\snl
We come to the following definition.

\inspr{2.4} Definition \rm
An {\it algebraic quantum group} is a regular multiplier Hopf
algebra with integrals.
\einspr

Algebraic quantum groups have been studied in [VD6].  The main result is 
duality:

\inspr{2.5} Theorem \rm
Let $(A,\Delta)$ be an algebraic quantum group.  Let $\varphi$
be a left integral.  Set $\hat A = \{ \varphi (\,\cdot\, a) \mid  a \in A \}$.  Considering
$\hat A$ as a subspace of $A^\prime$, the dual of the coproduct and the product
of $A$ give a product and a coproduct $\hat \Delta$ on $\hat A$ such that $(\hat
A, \hat \Delta)$ is again an algebraic quantum group.
\einspr

We call $(\hat A, \hat \Delta)$ the dual of $(A,\Delta)$.  It is also shown in
[VD6] that the dual of $(\hat A, \hat \Delta)$ is canonically isomorphic with
$(A,\Delta)$.
\snl
For algebraic quantum groups, we can push the result of 2.2 a little further:

\inspr{2.6} Proposition \rm
Let $(A,\Delta)$ be an algebraic quantum group. Given
elements $(a_1,...,a_n)$ in $A$, there exists an element $e \in A$ such that
$a_ie = ea_i = a_i$ for all $i$.
\snl\bf
Proof : \rm
Let $\omega_i, \rho_i$ be linear functionals so that
$$\sum \omega_i (ba_i) + \sum \rho_i (a_ib) = 0$$
for all $b \in A$.  Then, as in the proof of 2.2, we have to show now that
$$\sum \omega_i (a_i) + \sum \rho_i (a_i) = 0.$$
But, when $\sum \omega_i (ba_i) + \sum \rho_i (a_i b) = 0$ for all $b \in A$,
then also
$$\sum (\omega_i \otimes \iota)((\Delta (b)(a_i \otimes 1)) + \sum (\rho_i
\otimes \iota)((a_i \otimes 1) \Delta (b)) = 0$$
for all $b \in A$.  If we apply a left integral $\varphi$ we get $\varphi (b)
(\sum \omega_i (a_i) + \sum \rho_i (a_i)) = 0$.  So, just take any $b$ so that
$\varphi (b) \ne 0$. \hfill $\blacksquare$
\einspr
Observe that the proof here is simpler than in 2.2 and that the result is
slightly stronger.  We have to mention that 2.6 first has been proved by J.
Kustermans, but his proof was more complicated.
\snl
Again, this stronger result makes the use of the Sweedler notation for algebraic
quantum groups even more acceptable.  Consider e.g.\ an element $a$ and
a finite number of elements
$a_i,b_i,c_i,d_i$.  Then we can choose one element $e$ so that $ea_i = a_i$ for
all $i$ and
$$\align
\Delta (a) (b_i \otimes 1) & = (\Delta (a)(b_i \otimes 1))(1 \otimes e) \\
& = (\Delta (a)(1 \otimes e))(b_i \otimes 1) \\
(1 \otimes c_i)\Delta (a) & = ((1 \otimes c_i) \Delta (a)) (1 \otimes e) \\
& = (1 \otimes  c_i)(\Delta (a)(1 \otimes e))   \\
(d_i \otimes 1)\Delta (a) & = ((d_i \otimes 1) \Delta (a))(1 \otimes e) \\
& = (d_i \otimes 1)(\Delta (a)(1 \otimes e)).
\endalign
$$
for all $i$.  So, if we use $\sum a_{(1)} \otimes a_{(2)}$ for $\Delta (a)(1
\otimes e)$, we can write
$$\align
\sum a_{(1)} \otimes a_{(2)} a_i & = \Delta (a) (1 \otimes a_i) \\
\sum a_{(1)} b_i \otimes a_{(2)} & = \Delta (a) (b_i \otimes 1) \\
\sum a_{(1)} \otimes c_i a_{(2)} & = (1 \otimes c_i)\Delta (a) \\
\sum d_i a_{(1)} \otimes a_{(2)} & = (d_i \otimes 1) \Delta (a)
\endalign
$$
for all $i$.
\snl
In connection with the other remark we made in the general case, we can mention
that for algebraic quantum groups, the set of functionals $\omega \in A^\prime$
such that $(\omega \otimes \iota) \Delta (a)$ and $(\iota \otimes \omega)\Delta
(a)$ are in $A$ for all $a \in A$ precisely coincides with $M(\hat A)$ (see
[K]).
\snl
Algebraic quantum groups will play an important r\^ole in this paper. In
fact, the main duality result in section 7 will only be proven for an
algebraic quantum group.
\nl
For completeness, we now recall also the notions of
algebraic quantum groups of discrete and compact type.
These two special cases will not play a special r\^ole here.
Nevertheless, they are the best known types of algebraic quantum groups
and so, it is good to recall their definitions.
\snl
First recall the definition of a cointegral.

\inspr{2.7} Definition \rm
Let $(A,\Delta)$ be a multiplier Hopf algebra.  A non-zero
element $h \in A$ is called a {\it left cointegral} if $ah = \epsilon(a)h$ for
all $a \in A$.  Similarly, a non-zero element $k$ in $A$ is called a
{\it right cointegral} if $ka = \epsilon (a)k$ for all $a \in A$.
\einspr

It has been shown in [VD-Z1] that also cointegrals are unique (up to a scalar) for
regular multiplier Hopf algebras. They are faithful in the sense that $(\omega
\otimes \iota)\Delta (h) = 0$ and $\omega \in A^\prime$ implies $\omega = 0$ and
$(\iota \otimes \omega) \Delta (h) = 0$ implies $\omega = 0$.  Again the
antipode will convert a left cointegral into a right one and a right one into a
left one.
\snl
Regular multiplier Hopf algebras with cointegrals have been studied in
[VD-Z1].
There it is show that integrals exist in this case.  We also have obtained some
necessary and sufficient conditions on the underlying algebra $A$ for
$(A,\Delta)$ to have cointegrals.  We have the following definition and result.

\inspr{2.8} Definition \rm
A regular multiplier Hopf algebra with cointegrals is called
a multiplier Hopf algebra of discrete type (or an algebraic quantum group of
discrete type).  A Hopf algebra with integrals is called (an algebraic quantum
group) of compact type.
\einspr

It is shown in [VD6] that the dual of a discrete type is of compact type and of
course that the dual of a compact type is of discrete type.  Because for the
discrete type, we have integrals and cointegrals, we get a class of multiplier
Hopf algebras that shares most properties with the finite-dimensional Hopf
algebras.
\snl
In the discrete type case, we have shown in [VD-Z1] that the elements $e$ obtained
in proposition 2.6 can be chosen to be idempotent.
The existence of these {\it local units} is obvious in the case of a
discrete group $G$. Then, $A$ is the algebra $K(G)$ of complex functions
with finite support on $G$ (and pointwise operations). The local units
are the functions that are $1$ on a finite subset of $G$ and $0$ outside
this subset. Similarly in the case of a discrete quantum group in the
sense of [E-R] and [VD4]. In this case, the underlying algebra is a
direct sum of matrix algebras and the local units are the sums of
finitely many identities of these components. For another example, see
[VD6] and [VD-Z1].
\snl
In the discrete type case, we also have
that $M(\hat A) = \hat A$ as $\hat A$ has an identity.  So, the set of
functionals $\omega \in A^\prime$ satisfying $(\omega \otimes \iota)\Delta(a)
\in A$ and $(\iota \otimes \omega)\Delta(a) \in A$ for all $a \in A$ are
precisely those in $\hat A$.
\abs
\abs
\head 3. Multiplier Hopf Algebra Modules\endhead
\abs
In this section we will study some basic properties of left $A$-modules where
$A$ is the underlying algebra of a multiplier Hopf algebra.  In particular, in
what follows, $A$ will be an algebra over $\Bbb C$, with or without an identity
but with a non-degenerate product and satisfying property 2.2 of the previous
section.
\snl
We consider a left $A$-module $R$.  So, $R$ is a vector
space over $\Bbb C$ and
we have a bilinear map $(a,x) \in A \times R \rightarrow ax \in R$ satisfying
$(aa^\prime)x = a(a^\prime x)$ for all $a,a^\prime \in A$ and $x \in R$.
\snl
Left $A$-modules where already considered in [D-VD]. Here, we will use a
slightly different terminology.
\snl
If the algebra $A$ has an identity, it is normal to assume that $1 x = x$
for all $x \in R$.  This means that the module is unital.  In our context, it
seems to be most natural to extend this notion in the following way:

\inspr{3.1} Definition \rm
Let $R$ be a left $A$-module.  We call $R$ {\it
unital} if $AR = R$.
\einspr

It is clear that this notion coincides with the condition $1x = x$ for all $x$
when $A$ has an identity.
\snl
We also have the following easy consequence.

\inspr{3.2} Proposition \rm
Let $R$ be a unital left $A$-module.  If $x \in R$ and $ax
= 0$ for all $a \in A$, then $x=0$.

\snl\bf Proof : \rm
So, let $x \in R$ and assume that $ax = 0$ for all $a \in A$.  By
assumption, we can write $x = \sum a_i x_i$ where $a_i \in A$ and $x_i \in R$
for all $i$.  Choose $e \in A$ such that $ea_i = a_i$ for all $i$.  Then $ex =
x$.  But $ex = 0$ and hence $x = 0$.
\einspr

We say that a unital $A$-module is {\it non-degenerate}.
\snl
Remark that
in  [D-VD] we call a module non-degenerate if R is unital and satisfies the
property in proposition 3.2. This ambiguity has to do with the notion of
a non-degenerate $^*$-representation $\pi$ of (say) a C$^*$-algebra $B$ on a Hilbert
space $\Cal H$. In that case, the density of $\pi(B)\Cal H$ in $\Cal H$
is equivalent with the non-degeneracy condition $\pi(B)\xi={0}$ implies
$\xi=0$ (which is essentially the condition in 3.2). See [P].
\snl
This observation might suggest that the converse of
proposition 3.2 is true. This is not the case. Consider $A^\prime$ as a left $A$
module by $a \omega = \omega (\,\cdot\, a)$.  This module is non-degenerate because
if $\omega (\,\cdot\, a) = 0$ for all $a$, then $\omega = 0$.  On the other hand, in
general, we do not have that $A^\prime = AA^\prime$.
\snl
So we see that in this purely algebraic context, the notion of
non-degeneracy (as formulated in proposition 3.2) is weaker than
the notion unital.
This is different in the topological context when dealing with
$^*$-algebras.
\snl
We believe that many results in this paper can still be obtained under
the weaker condition. In some cases however, it may be necessary to have
the stronger one. Also, for other reasons, it seems very natural to
assume the condition $AR=R$. We have e.g.\ that $AA=A$ and conditions
like $\Delta(A)(1\ot A)=A\ot A$. Also, when we consider a homomorphism
from one algebra into the multiplier algebra of another one, we will
require a similar condition (see definition 4.12 in the next section).
So, because it seems very natural and for convenience, we will assume
throughout the paper that we are dealing with unital modules. There are
a few exceptions but we will indicate these very clearly.
\snl
Remark that for a unital $A$-module $R$, for any $x \in R$ there is an
$e \in A$
such that $ex = x$.  This was used in the proof of the proposition.  In fact,
using the same argument, we have for all $x_1,x_2,...,x_n \in R$ an element $e$
such that $ex_i = x_i$ for all $i$. This is important for working with
the Sweedler notation: it means that elements in a unital left
$A$-module will cover elements $a_{(k)}$.
\snl
Next we show that unital $A$-modules can be extended to modules over $M(A)$.
This result can already be found in [D-VD].

\inspr{3.3} Proposition \rm
Let $R$ be a unital left $A$-module.  Then, there is a
unique extension to a left $M(A)$-module.  We have $1x = x$ for all $x \in R$
where now $1 \in M(A)$.

\snl\bf Proof : \rm
If $M(A)$ also acts on $R$ and if this action extends the action of $A$,
we must have $m(ax) = (ma)x$ for all $x \in R$, $a \in A$ and $m \in M(A)$.
Because $R$ is unital, this relation completely determines the action of $M(A)$
and we must have $1x = x$ for all $x$.  On the other hand, we can use this
formula to define the action of $M(A)$.  To prove that this is well-defined,
assume that $a_i \in A$ and $x_i \in R$ and that $\sum a_i x_i = 0$.  Choose $e
\in A$ so that $ea_i = a_i$ for all $i$.  For any $m \in M(A)$ we have
$$\sum (m a_i) x_i  = \sum (me) (a_i x_i)  = (me) \sum a_i x_i = 0. $$
It follows that we can define the action of $M(A)$ by $m(ax) = (ma)x$ whenever
$a \in A$, $x \in R$ and $m \in M(A)$. \hfill $\blacksquare$
\einspr

Suppose that we did not have an algebra with {\it local units} (as in
propositioin 2.2).  Then also proposition 3.3 might fail.  However, if $R$ is
both unital and non-degenerate, we still have 3.3.  Indeed, observe that, if
$\sum a_i x_i = 0$, then for all $b \in A$ and $m \in M(A)$ we get
$$b(\sum (ma_i)x_i)= \sum (bma_i) x_i = (bm) \sum a_i x_i = 0\,.$$
Because the module is non-degenerate we obtain $\sum (ma_i) x_i = 0$ 
(see again [D-VD]). We will come back to this problem in
section 5 where we
will also consider modules over the smash product (which need not to
have such local units).
\snl

Now, let $(A,\Delta)$ be a regular multiplier Hopf algebra.  Then we can use
$\Delta$ to define tensor products of unital left $A$-modules.  Indeed,
suppose that $R$ and $T$ are unital left $A$-modules. Then we can make
$R \otimes T$ into a left $(A \otimes A)$-module by $(a \otimes b)(x \otimes y) =
ax \otimes by$.  It is clear that also this module will be unitial.  Hence, we
can use proposition 3.3 to extend it to $M(A \otimes A)$.  Then it is possible
to use $\Delta$ to get an action of $A$ on $R \otimes T$ given by
$$a(x \otimes y) = \Delta (a)(x \otimes y).$$
Now, because $\Delta (A) (A \otimes A) = A \otimes A$ we get
$$\align
\Delta (A) (R \otimes T) & = \Delta (A) (A \otimes A) (R \otimes T) \\
& = (A \otimes A)(R \otimes T) \\
& = R \otimes T
\endalign
$$
so that also the action of $A$ on $R \otimes T$ is a unital action.
\snl
In fact, we get the following result.

\inspr{3.4} Proposition \rm
Let $A$ be a regular multiplier Hopf algebra.  Denote by ${\Cal
M}$ the category of unital left $A$-modules and morphisms.  Then ${\Cal M}$ is a
monoidal category with unit (for the product structure defined above and with
the diagonal $A$-module structure).
\einspr

The unit is of course given by $\Bbb C$. The module structure in this
case is $a\lambda=\epsilon(a)\lambda$ whenever $a\in A$ and $\lambda\in
\Bbb C$. It is clear that this is a unit because $\Delta(a)(\lambda\ot
x)=\lambda(ax)=a(\lambda x)$ for all $a\in A$, $x\in R$ and $\lambda \in
\Bbb C$.
\nl
In prinicple, when $R$ is a unital left $A$-module, we can also define the
{\it fixed points} in $R$ as the elements $x\in R$ that satisfy
$ax=\epsilon(a)x$ for all $a\in A$. This notion is a little too
restrictive. Think e.g.\ of the algebra of complex functions with finite
support on a group $G$ and let it act on itself by left convolution. This
will give a unital action of the group Hopf algebra $\Bbb C$. There will
be no fixed points. Only constant functions would be fixed points, but
these are not in our algebra. So, in fact, it is more natural to look at
fixed points in the multiplier algebra of $R$. But to do this, one first
has to extend the action of $A$ to this multiplier algebra. This will be
done in section 4.
\snl
For the module actions, we will use different notations. Here, we have
used $ax$ for the action of $a$ on $x$. In some papers, $a.x$ is used
instead. Sometimes, certainly when more different types of
actions are involved, we will also use $a\triangleright x$.
\abs
\abs
\head 4. Module Algebras over Multiplier Hopf Algebras\endhead
\abs
Again, let $(A, \Delta)$ be a regular multiplier Hopf algebra.  Now assume that
$R$ is an algebra over $\Bbb C$, with or without identity but with non-degenerate
product.  Assume that $R$ is a unital left $A$-module, cf.\ definition
3.1.
\nl
The main definition in this section is the following.

\inspr{4.1} Definition \rm
We say that $R$ is a left $A$-module algebra if
$$a(xx^\prime) = \sum (a_{(1)} x)(a_{(2)}x^\prime)$$
for all $a \in A$ and $x,x^\prime \in R$.
\einspr

Let us first discuss the right hand side of the above equality.  We have
seen in the previous section that the elements $x$ and $x'$ can be used
to cover $a_{(1)}$ and $a_{(2)}$. But it is also possible to view this
expression as
$m(\Delta (a)(x \otimes x^\prime))$
where now $m$ denotes multiplication in $R$.  Indeed, we have seen in
proposition 3.3 that this makes sense because $R$ is unital.  Then, the
condition in 4.1 can be rewritten as $m\Delta(a) = am$, saying that $m$
is an $A$-module map.
\snl
When $R$ is a left $A$-module algebra, we will also say that the
multiplier Hopf algebra $A$ {\it acts on} $R$. This last formulation is
closer to the notion of an action of a group (which in fact is the
original situation). See example 4.3 below.
\snl
Before we give some examples however,
let us first prove some immediate consequences
of the definition. In fact, in the following lemma, we give some
formulas that are essentially equivalent with the formula in the
definition above.

\inspr{4.2} Lemma \rm
For any $a\in A$ and $x,x'\in R$ we have
$$\align (ax)x'&= \sum a_{(1)} (x (S(a_{(2)})x')) \\
         x(ax')&= \sum a_{(2)} ((S^{-1}(a_{(1)})x)x').
  \endalign$$

\snl\bf Proof : \rm
We have
$$\align \sum a_{(1)} (x (S(a_{(2)})x^\prime))
         &=\sum (a_{(1)}x) (a_{(2)}S(a_{(3)})x') \\
         &=\sum \epsilon(a_{(2)})(a_{(1)}x)x'=(ax)x'
   \endalign$$
and similarly for the other formula.  \hfill $\blacksquare$
\einspr

Remark that $x'$ covers $a_{(2)}$ in the first formula and that
$a_{(1)}$ is covered by $x$ in the second formula of the lemma.
\snl
We will use these formulas at several places in the paper. In
particular, we will use them later in this section to extend the action
of $A$ on $R$ to $M(R)$.
\snl
Let us now consider some examples.
We start with the natural one.
There is nothing new about this example
because we have an example of an action of a Hopf algebra. But it serves
as a motivation, also for actions of multiplier Hopf algebras.

\inspr{4.3} Example \rm
Let $G$ be a group and let $A$ be the group algebra $\Bbb CG$.
Denote the imbedding of $G$ in $A$ by $p \rightarrow \lambda_p$ so that $A$ is
spanned by the elements $\{ \lambda_p \mid p \in G \}$ and that $\Delta$ is
given by $\Delta (\lambda_p) = \lambda_p \otimes \lambda_p$ for all $p \in G$.
Now, let $\alpha$ be an action of $G$ on an algebra $R$ by means of automorphisms
of $R$.  We associate an action of $A$ on $R$ by $\lambda_p x = \alpha_p(x)$
whenever $p \in G$ and $x \in R$.  Because $\alpha$ is an action of $G$, we have
that $R$ is an $A$-module.  Since also $\alpha_e(x) = x$ for all $x \in R$
when $e$ is the identity in $G$, we have a unital $A$-module.  Finally, because
$\alpha_p$ is an automorphism, we get
$$\align
\lambda_p (xx^\prime) & = \alpha_p (xx^\prime) = \alpha_p (x) \alpha_p
(x^\prime) \\
& = (\lambda_p x) (\lambda_p x^\prime).
\endalign
$$
So, $R$ is a $A$-module algebra.
\einspr
The next example is (in general) not a Hopf algebra example.

\inspr{4.4} Example \rm
Let $G$ be group and now let $A$ be the algebra of complex
functions with finite support in $G$.  Let $\Delta$ be defined by $(\Delta
f)(p,q) = f(pq)$ for all $p,q \in G$ and $f \in A$.  Now denote by $\delta_p$
the function that is $1$ on $p$ and $0$ everywhere else.  Suppose  that we
have an action of $A$ on an algebra $R$.  Denote $R_p = \delta_p R$.  It is
clear that $\delta_p$ acts as a projection map from $R$ to the space $R_p$.
Moreover, if $R$ is unital, we will get that $R$ can be identified with the
direct sum of these subspaces $R_p$ with $p \in G$.  Because
$$\Delta (\delta_p) = \sum_{p=qr} \delta_q \otimes \delta_r$$
we have that
$$\delta_p(xx^\prime) = \sum_{p=qr} (\delta_q x)(\delta_r x^\prime).$$
In particular, if $x \in R_q$ and $x^\prime \in R_r$ we get $xx^\prime \in
R_{qr}$.  So $R$ is a $G$-graded algebra.  In fact, also conversely, if $R$ is a
$G$-graded algebra, with grading $R = \oplus_{p \in G} R_p$, then $A$ acts on
$R$ if we let $\delta_p x = x$ when $x \in R_p$ and $\delta_p x = 0$ when $x \in
R_q$ and $q \ne p$.
\einspr

Let us now consider the adjoint action. It will play a fundamental
r\^ole in what follows.

\inspr{4.5} Proposition \rm
Let $A$ act on itself by the map $a \otimes x
\rightarrow \sum a_{(1)} x S(a_{(2)})$ from $A\ot A$ to $A$.
This action makes $A$ into an $A$-module algebra.
\snl\bf Proof : \rm 
First observe that the above expression makes sense.
Also, compare this formula with the first one
in lemma 4.2. Because for all $a,b,x \in A$ we get
$$\align
\sum (ab)_{(1)} xS((ab)_{(2)})
& = \sum a_{(1)} b_{(1)} xS(a_{(2)} b_{(2)}) \\
& = \sum a_{(1)} b_{(1)} xS(b_{(2)}) S(a_{(2)}),
\endalign
$$
we clearly have a left $A$-module structure. To see that we have a unital
module, observe that the surjectivity of the map
$$a \otimes x \rightarrow \Delta (a) (x \otimes 1)$$
will give that $pS(q)$ is in the range of this module map for all
$p,q\in A$ and so is $A$.
Finally, for all $a,x,x^\prime \in A$, we also get
$$\align
\sum (a_{(1)} x S(a_{(2)}))(a_{(3)} x^\prime S(a_{(4)})) & = \sum
\epsilon (a_{(2)}) a_{(1)} xx^\prime S(a_{(3)}) \\
& = \sum a_{(1)} xx^\prime S(a_{(2)})
\endalign
$$
so that we have an $A$-module algebra structure.  \hfill $\blacksquare$
\einspr

Of course, if we consider the adjoint action of proposition 4.5 for the
group algebra (example 4.3), we
get the adjoint action of $G$ on $\Bbb CG$ given by $\alpha_p(x) = \lambda_p x
\lambda^{-1}_p$.
\snl
Also the following example is closely related with the adjoint action.

\inspr{4.6} Example \rm
Start with a left
$A$-module $R$. Let $L$ be the vector space of linear maps from $R$ to
$R$. Then $L$ is an algebra with identity (the product is the
composition of maps). The action of $A$ yields an action on $L$ given by
the following formula. If $\lambda\in L$, let $a\lambda$ be the linear
map on $R$ defined by
$$(a\lambda)x=\sum a_{(1)}\lambda (S(a_{(2)})x).$$
Remark that $x$ will cover $a_{(2)}$. Also here, compare this formula
with the first formula in lemma 4.2. Similarly as in example 4.5,
we have made $L$ into an $A$-module algebra.
\einspr

Now we want to extend the action of $A$ to the multiplier algebra
$M(R)$. Since elements in the multiplier algebra $M(R)$ yield linear
maps from $R$ to $R$, we can use the formula in example 4.6 to extend
the action of $A$ from $R$ to $M(R)$. This is done in the following
proposition. Remark that again, we are motivated by the formulas in
lemma 4.2.

\inspr{4.7} Proposition \rm
We can define  $am \in M(R)$ for all $a \in A$ and $m \in
M(R)$ by
$$\align
(am)x & = \sum a_{(1)} (m(S(a_{(2)})x)) \\
x (am) & = \sum a_{(2)} ((S^{-1} (a_{(1)})x)m)
\endalign
$$
whenever $x \in R$.
\snl\bf Proof : \rm
We will check the relation $(x^\prime (am))x^{\prime \prime} = x^\prime
((am)x^{\prime \prime})$ with the expressions in the formulation of the
proposition. For the left hand side, we get
$$(x^\prime (am))x^{\prime \prime} = \sum (a_{(2)} ((S^{-1} (a_{(1)}
)x^\prime)m))x^{\prime \prime}.$$
For the right hand side we have
$$\align
x^\prime ((am)x^{\prime \prime}) & = \sum x^\prime (a_{(1)}
(m(S(a_{(2)})x^{\prime \prime}))) \\
& = \sum (a_{(2)} S^{-1}(a_{(1)})x^\prime)(a_{(3)} (m(S(a_{(4)})x^{\prime
\prime}))) \\
& = \sum a_{(2)} ((S^{-1}(a_{(1)})x^\prime) (m(S(a_{(3)})x^{\prime\prime}))) \\
& = \sum (a_{(2)} ((S^{-1}(a_{(1)})x^\prime)m))(a_{(3)} S(a_{(4)})x^{\prime
\prime}) \\
& = \sum (a_{(2)} ((S^{-1}(a_{(1)} )x^\prime)m)) x^{\prime \prime}
\endalign
$$
And we find the same. \hfill$\blacksquare$
\einspr

Remark that in all the formulas, appearing in the proof above, enough elements
$a_{(k)}$ are covered.  In particular, in the defining relations e.g.\ $a_{(2)}$
is covered by $x$ in the first one and $a_{(1)}$ by $x$ in the
second one.
\snl
We have some obvious results on this extended action.

\inspr{4.8}  Proposition \rm
We have $(aa^\prime)m = a(a^\prime m)$ for all $a,a^\prime
\in A$ and $m \in M(R)$. Also $a1 = \epsilon (a)1$.
\snl\bf Proof : \rm
Take $a,a^\prime \in A$ and $m \in M(R)$.  For all $x \in R$ we get
$$\align
((aa^\prime)m) x & = \sum (aa^\prime)_{(1)} (m(S((aa^\prime)_{(2)})x)) \\
& = \sum a_{(1)} a^\prime_{(1)} (m(S(a^\prime_{(2)}) S(a_{(2)})x)) \\
& = \sum a_{(1)} ((a^\prime m)(S(a_{(2)}) x)) \\
& = (a(a^\prime m))x.
\endalign
$$
This proves the first statement.
\snl
To prove the second statement, take $a\in A$ and $x\in R$. Then
$$(a1)x=\sum a_{(1)}(S(a_{(2)})x)=\sum(a_{(1)}S(a_{(2)})x=\epsilon(a)x.$$
Hence $a1=\epsilon(a)1$.\hfill $\blacksquare$
\einspr

So, we obtain that $M(R)$ is a left $A$-module.
One might guess
that still $M(R)$ is a left $A$-module algebra. But there is a problem.
One can no longer expect that $M(R)$ is a unital left $A$-module. In
general, the condition $AM(R)=M(R)$ will be much too strong. So, we can
not give a meaning to the formula
$$a(mm')=\sum(a_{(1)}m)(a_{(2)}m')$$
for all $a\in A$ and $m,m'\in M(R)$ in the same way as we did for $R$.
\snl
Fortunately, we do have that the action of $A$ on $M(R)$ is still
non-degenerate:

\inspr{4.9} Proposition \rm
If $m\in M(R)$ and $am=0$ for all $a$, then $m=0$.
\snl\bf Proof : \rm
For all $a\in A$ and $x\in R$, we have
$$\align
m(S(a)x) &=\sum S(a_{(1)})a_{(2)} (m(S(a_{(3)})x)\\
          &=\sum S(a_{(1)})((a_{(2)} m)x)=0.
\endalign$$
As this is true for all $a\in A$ and $x\in R$, and because $R$ is
unital, we have $m=0$.
\hfill $\blacksquare$
\einspr

Now, we are ready to study the fixed points.
\snl
We have already mentioned in section 3 that fixed points can be defined
in $R$ as the elements $x$ satisfying $ax=\epsilon(a)x$ for all $a\in
A$. They obviously form a subalgebra of $R$. But the set of fixed points
in $R$ may be too small in general.
\snl
It makes more sense to look at fixed points in the multiplier algebra
$M(R)$.
This is possible since we have extended the action of $A$ to $M(R)$.
This gives the
following definition.

\iinspr{4.10} Definition \rm
A {\it fixed point} in $M(R)$ is an element $m$ satisfying $am=\epsilon(a)m$
for all $a\in A$.
\einspr

It is not so hard to prove the following property.

\iinspr{4.11} Proposition \rm
If $m$ is a fixed point in $M(R)$, then
$$\align  a(mx)&=m(ax)\\
          a(xm)&=(ax)m
  \endalign$$
for all $a\in A$ and $x\in R$.

\snl\bf Proof : \rm
Because $m$ is a fixed point, we get for all $a,a'\in A$ and $x\in R$
$$\align a'S(a)(mx)&=\sum a'S(a_{(1)})(a_{(2)}m)x\\
                 &=\sum a'S(a_{(1)})a_{(2)}(m(S(a_{(3)})x))\\
                 &=a'm(S(a)x).
  \endalign$$
We use $a'$ for covering purposes. Now we can cancel $a'$ and replace
$S(a)$ by $a$ to get $a(mx)=m(ax)$ for all $a\in A$ and $x\in R$. This proves the first
formula.
Similarly
$$\align a'S^{-1}(a)(xm)&=\sum a'S^{-1}(a_{(2)})(x(a_{(1)}m))\\
         &=\sum a'S^{-1}(a_{(3)})a_{(2)}((S^{-1}(a_{(1)})x)m)\\
         &=a'(S^{-1}(a)x)m
  \endalign$$
and we get the second formula.
\hfill $\blacksquare$
\einspr

This means that multiplication with $m$, both from the left and from the
right, commutes with the action of $A$ on $R$. It is not so hard to see
that these facts imply that $m$ is a fixed point. So, we could also use
this property to define fixed points in $M(R)$ and we do not really need
the extension of the action to introduce this notion. It also follows
that the fixed points in $M(R)$ form actually a subalgebra. This is not
completely obvious from the way we defined it as we do not have that
$M(R)$ is an $A$-module algebra.
\snl
We can also look for fixed points in the algebra $L$ of linear maps from
$R$ to $R$ for the action considered in example 4.6. So, a linear map
$\lambda$ will be a fixed point if
$$\sum a_{(1)}\lambda(S(a_{(2)})x)=\epsilon(a)\lambda (x)$$
and just as in the proof of the previous proposition, we get that this
is equivalent with $a\lambda(x)=\lambda(ax)$ for all $a$ and for all
$x$.
\nl
We now want to introduce the notion of an inner action. To motivate the
definition (definition 4.13 below) let us again first consider an action
$\alpha$ of a group $G$ on an algebra $R$ with identity.  Then $\alpha$ is
called inner if there is a homomorphism $p \rightarrow u_p$ from $G$ into the
invertible elements of $R$ so that $\alpha_p(x) = u_p xu^{-1}_p$ for all $p \in
G$ and $x \in R$.  It is no strong restriction to require $u_e = 1$
(because then $\alpha_e = 1$).  If $R$ has no identity, this notion would
be too restrictive and it is much more natural to allow homomorphisms from $G$
to $M(R)$.
\snl
In the case of a regular multiplier Hopf algebra, we will have to consider
therefore homomorphism $\gamma$ of $A$ into $M(R)$.  But also here we have to
require them to be unital in the following sense.

\iinspr{4.12} Definition \rm
A homomorphism $\gamma :  A \rightarrow M(R)$ is called
unital if $\gamma(A)R = R\gamma(A) = R$.
\einspr

In [VD3] and [D-VD], we spoke about non-degenerate homomorphisms. Again,
for this, we were inspired by the C$^*$-algebra context. But just as for
modules, it now seems to be more natural to call this property unital.
Indeed, just as for actions, unital homomorphisms have unique extensions to
homomorphisms
from $M(A)$ to $M(R)$.  This extension is still denoted by $\gamma$ and it can
be defined by e.g.\ $\gamma(m)(\gamma(a)x) = \gamma(ma)x$ and $(x \gamma(a)) \gamma
(m) = x \gamma (am)$ whenever $a \in A$, $m \in M(A)$ and $x \in R$.
For this extension we have $\gamma(1)  = 1$.  That is why, also here, we call
$\gamma$ unital.
\snl
Now, we are ready to give our definition of an {\it inner action}.

\iinspr{4.13} Definition \rm
Let $(A, \Delta)$ be a regular multiplier Hopf algebra acting
on an algebra $R$.  We call the action inner if there is a unital homomorphism
$\gamma : A \rightarrow M(R)$ so that
$$ax = \sum \gamma (a_{(1)}) x \gamma (S(a_{(2)}))$$
for all $a$ in $A$ and $x \in R$.
\einspr

Remark that the expression on the right makes sense because, as $\gamma$ is
unital, for all $x \in R$ we have $e \in A$ so that $x = \gamma (e)x$.  This
element $e$ covers $a_{(1)}$.
\snl
It is obvious that the adjoint action, as defined in proposition 4.5, is
actually inner (with $\gamma = \iota)$.
\snl
A natural extension of innerness is the cocycle equivalence of two
actions (sometimes called weak equivalence). This type
of equivalence was first considered (and appeared very naturally) in the
context of actions of locally compact groups on von Neumann algebras
(see e.g.\ [VD1]). It is important in crossed product theory.
\snl
To motivate the definition (definition 4.14 below), let us again first
have a look at two actions $\alpha$ and $\beta$ of a group on an
algebra $R$ with identity.  They are called cocycle (or weak)
equivalent if there is a map
$p \rightarrow u_p$ of $G$ into the invertible elements of $R$ satisfying
\snl
\quad i) $u_{pq} = u_p \alpha_p (u_q)$ for all $p,q \in G$,
\snl
\quad ii) $\beta_p (x) = u_p \alpha_p (x) u^{-1}_p$ for all  $p \in G$
and $x \in R$.
\snl
Again, it is natural to require $u_e = 1$.
Remark that condition i) is also very natural when $\alpha$ and $\beta$ are
actions satisfying ii).  Indeed, for $p,q \in G$ we get
$$\align
\beta_{pq} (x) & = u_{pq} \alpha_{pq} (x) u^{-1}_{pq} \\
& = u_p \alpha_p (u_q) \alpha_p (\alpha_q (x)) \alpha_p (u_q)^{-1} u^{-1}_p \\
& = u_p \alpha_p (u_q \alpha_q (x) u^{-1}_q) u^{-1}_p \\
& = \beta_p (\beta_q (x)).
\endalign
$$
In order to formulate this notion in the context of multiplier Hopf algebras,
we will have to replace the map $p \rightarrow u_p$ by a linear map from
$A$ to $M(R)$. We are able to do this since we have already defined the
extension of the action of $A$ to $M(R)$. We only give the definition in
the case of a Hopf algebra (for reasons we will explain afterwards).

\iinspr{4.14} Definition \rm
Let $(A,\Delta)$ be a Hopf algebra, acting on an
algebra $R$ in two ways. Denote these two actions by $a\triangleright_1
x$ and $a\triangleright_2 x$. We call these two actions cocycle
equivalent if there is a linear map $\gamma : A\to M(R)$ such that
$\gamma(1)=1$ and satisfying
\snl
\quad i) $\gamma(aa')=\sum \gamma (a_{(1})(a_{(2)}\triangleright_1
\gamma(a'))$
\snl
\quad ii) $\sum (a_{(1)}\triangleright_2 x )\gamma(a_{(2)})
=\sum\gamma(a_{(1)})(a_{(2)}\triangleright_1 x)$
\snl
for all $a,a'\in A$ and $x\in R$.
\einspr

There are some difficulties to extend this notion to the case of regular
multiplier Hopf algebras. We will need some kind of non-degeneracy
condition on $\gamma$ (to replace the condition $\gamma(1)=1$),
but it is not so clear how to define this since
$\gamma$ is no longer assumed to be a homomorphism. It is however
possible to give a meaning to the two other reguirements. The second one
is easy as the actions are assumed to be unital. To give a meaning to
the first condition, we can rewrite it as
$$
\gamma(aa')x=\sum \gamma
(a_{(1})a_{(2)}\triangleright_1((\gamma(a')(S(a_{(3)})\triangleright_1
x)))
$$
for all $a,a'\in A$ and $x\in R$.
\FF
\head 5. Smash Products\endhead
\abs
Let $(A,\Delta)$ be a regular multiplier Hopf algebra and let it act on an
algebra $R$ with a non-degenerate product as in the previous section.
\snl
By assumption, $R$ is a left $A$-module and we have
$$a(xy) = \sum (a_{(1)}x)(a_{(2)}y)$$
for all $x,y \in R$ and $a \in A$.  If we also think of $R$ as a left R-module,
we are led to the following definition.

\inspr{5.1} Definition \rm
Let $V$ be a vector space which is both a left $A$-module
and a left $R$-module.  We say that $V$ is a covariant $A$-$R$-module if
$$a(xv) = \sum (a_{(1)}x)a_{(2)}v$$
for all $a \in A$, $x \in R$ and $v \in V$.
\einspr

Remark that this formula makes sense because $a_{(1)}$
is covered by $x$.  If $V$ is a unital left $A$-module, we can rewrite this
formula as
$$(ax)v = \sum a_{(1)} (x S(a_{(2)})v).$$
In this case, $a_{(2)}$ is covered by $v$. This reminds us of the
formulas for the adjoint action (see example 4.5 and 4.6).
These formulas are indeed very natural.
If we consider e.g.\ $R$ not only as a
left $A$-module, but also as a left $R$-module, then $R$ becomes a
covariant $A$-$R$-module and the two formulas above are found in the
definition 4.1 and lemma 4.2.
\snl
The name comes from the theory of actions of locally compact groups on
operator algebras (see e.g.\ [P]). The following examples refers to this
original notion of covariant representation (in the simpler setting of a
group action).

\inspr{5.2} Example \rm
Let $G$ be a group and assume that $\alpha$ is an action of $G$
on $R$.  Let $V$ be a vector space over $\Bbb C$.  Suppose that $\pi$ is a
representation of $R$ on $V$ and that $p \rightarrow u_p$ is a representation of
$G$ on $V$ satisfying
$$\pi(\alpha_p(x)) = u_p \pi (x)u^{-1}_p$$
for all $x \in R$ and $p \in G$.  Such a pair is called a covariant
representation of the covariant system $(R,G,\alpha)$ in [P].  Now let $A$ be
the group algebra $\Bbb CG$ of $G$.  Then, the respresentation of $G$ makes
$V$ into  a $A$-module.  Similarly, the representation of $R$ makes $V$ into a
$R$-module.  And the relation $\pi (\alpha_p(x)) = u_p \pi (x) u^{-1}_p$
precisely means that $V$ is a covariant $A$-$R$-module.
\einspr

Now, let $V$ be a covariant $A$-$R$-module.  It follows from the definition that
$$x ax'a'v = \sum (x(a_{(1)} x'))(a_{(2)} a' v)$$
when $a,a'\in A$, $x,x'\in R$ and $v\in V$.
This formula suggests a multiplication on the tensor product $R \otimes A$
of the spaces $R$ and $A$ making $V$ into a left $(R \otimes A)$-module.
\snl
This motivates the introduction of the smash product.

\inspr{5.3} Definition \rm
Define a multiplication on $R \otimes A$ by
$$(x \otimes a)(x' \otimes a') = \sum x(a_{(1)}x') \otimes a_{(2)} a'$$
whenever $x,x' \in R$ and $a,a' \in A$.
\einspr

There is no problem with this definition: $a_{(1)}$ is covered by $x'$ and
$a_{(2)}$ is covered by $a'$.
\snl
The product above is given by the map
$$(m\ot m) (\iota\ot \Gamma \ot\iota)$$
where $\Gamma:A\ot R \to R\ot A$ is defined by
$$\Gamma(a\ot x)=\sum a_{(1)}x\ot a_{(2)}$$
and where $m$ is used to denote both the multiplication in $A$ and in
$R$, considered as linear maps from $A\ot A$ to $A$ and $R\ot R$ to $R$
respectively. The associativity of the product, as well as other
properties, can be obtained from properties of $\Gamma$. This method was
followed in [D-VD] and before in [VD-VK]. Here we will follow the more
direct approach. The calculations one has to make are esssentially the
same in both cases. The direct approach has the advantage of being more
self-contained. On the other hand, we will encounter the twist map
$\Gamma$ later and so, it is good to have in mind the above formula for
the product.
\snl
We first have to prove the following.

\inspr{5.4} Lemma \rm
The product defined above is associative.

\snl\bf Proof : \rm
Take $x,x^\prime, x^{\prime \prime} \in R$ and $a,a^\prime,a^{\prime
\prime} \in A$.  Then
$$\align
((x \otimes a)(x^\prime \otimes a^\prime))(x^{\prime \prime} \otimes a^{\prime
\prime}) & = \sum ((x(a_{(1)} x^\prime)) \otimes (a_{(2)} a^\prime))(x^{\prime
\prime} \otimes a^{\prime \prime}) \\
& = \sum x(a_{(1)} x^\prime)(a_{(2)} a^\prime_{(1)} x^{\prime \prime}) \otimes
a_{(3)} a_{(2)}^\prime a^{\prime \prime} \\
& = \sum x(a_{(1)} (x^\prime (a^\prime_{(1)} x^{\prime \prime})))
\otimes a_{(2)} a^\prime_{(2)} a^{\prime \prime} \\
& = (x \otimes a) \sum x^\prime (a^\prime_{(1)} x^{\prime \prime}) \otimes
a^\prime_{(2)} a^{\prime \prime} \\
& = (x \otimes a)((x^\prime \otimes a^\prime)(x^{\prime \prime} \otimes
a^{\prime \prime})).
\endalign
$$
\hfill $\blacksquare$
\einspr

So, we have made $R \otimes A$ into an associative algebra.  In this context, we
will use the following notations.

\inspr{5.5} Notation \rm
When we consider $R \otimes A$ with the above product, we will
use the notation $R \# A$ instead and elements $x \otimes a$ will be denoted by
$x \# a$.
\einspr

From the motivation above, we see that any covariant $A$-$R$-module $V$
gives rise to a $R\# A$-module structure on $V$.
We will come back to this in proposition 5.11.
In particular, this is
the case for the covariant $A$-$R$-module $R$. The action of $R\# A$ on
$R$ is given by $(x\# a)x'=x(ax')$ whenever $a\in A$ and $x,x'\in R$. If
this module was faithful, we could have used it to define $R\# A$.
Unfortunately, in general, this will not be the case. We will however
use this module often and in some special cases, it will be faithful.
\snl
It is easy to see that, in the case of the example 5.2, the algebra $R \# A$ is
the algebra of functions from $G$ to $R$ with finite support with the twisted
convolution product
$$(\xi\eta)(p) = \sum_q \xi(q) \alpha_q (\eta(q^{-1} p)).$$
So, $R \# A$ is the crossed product $R \times_\alpha G$ in the sense of
[P].
\snl
We will now investigate the algebra $R \# A$ further.  The first
(important) step is the following.

\inspr{5.6} Lemma \rm
The product in $R \# A$ is non-degenerate.

\snl\bf Proof : \rm
Suppose that $\sum x_i \# a_i \in R \# A$ and that $(\sum x_i \# a_i)(x
\# a) = 0$ for all $x \in R$ and $a \in A$.  Then, using the definition of the
product in $R \# A$ and the non-degeneracy of the product in $A$, we get
$$\sum x_i(a_{i(1)} x) \# a_{i(2)} = 0$$
for all $x \in R$.  Now, we can proceed as before (see the technique
used in the proof of proposition 2.2 or in the proof of proposition 3.4
in [VD6]).  Apply $\Delta$ and $S$,
multiply with $a,a^\prime$ and replace $x$ by $a^{\prime \prime}x$ to obtain
$$\sum x_i (a_{i(1)} a^{\prime \prime} x) \otimes S(a_{i(2)} )a \otimes a_{i(3)}
a^\prime = 0$$
for all $a,a^\prime,a^{\prime \prime} \in A$ and $x \in R$.  Now, we can
{\it replace} $a^{\prime \prime}$ by $S(a_{i(2)})a$ as these elements are forced (by
the use of $a$ and $a^\prime$) to lie in a finite-dimensional space.  This gives
$$\sum x_i (a_{i(1)} S(a_{i(2)}) ax) \otimes a_{i(3)} a^\prime = 0$$
and hence
$$\sum x_i (ax) \otimes a_i a^\prime = 0$$
for all $a,a^\prime \in A$ and $x \in R$.  We can again cancel $a^\prime$.  And
because $R$ is unital so that $AR = R$, we can also cancel $ax$ and we obtain
$\sum x_i \otimes a_i = 0$.
\snl
Now suppose that $(x \# a)(\sum x_i \# a_i) = 0$ for all $x \in R$ and $a \in
A$.  This means that
$$\sum x(a_{(1)} x_i) \otimes a_{(2)} a_i = 0$$
for all $x$ and all $a$.  If we multiply the second factor from the left with
$a'$ and use the fact that $(1 \otimes A)\Delta (A) = A \otimes A$, we get that
also
$$\sum x(ax_i) \otimes a' a_i = 0$$
for all $x \in R$ and $a,a' \in A$.  We can cancel $x$, cancel $a'$ and use the
non-degeneracy of the module $R$ to get again $\sum x_i \otimes a_i = 0$.
\hfill$\blacksquare$
\einspr

If $R$ has an identity, then we get
$$\align
(1 \# a)(1 \# a^\prime) & = \sum (a_{(1)} 1) \# a_{(2)} a^\prime \\
& = \sum \epsilon (a_{(1)}) 1 \# a_{(2)} a^\prime \\
& = 1 \# aa^\prime
\endalign
$$
and it follows that $a \rightarrow  1 \# a$ is a homomorphism of $A$ into $R \#
A$.  Similarly, if $A$ has an identity we get
$$(x \# 1)(x^\prime \# 1) = x(1  x^\prime) \# 1 = xx^\prime \# 1$$
and $x \rightarrow x \# 1$ is a homomorphism of $R$ into $R \# A$.  If $R$ and
$A$ have identities, we get moreover that
$$\align
(x \# 1)(1 \# a) & = x \# a \\
(1 \# a)(x \# 1) & = \sum a_{(1)} x \# a_{(2)}
\endalign
$$
for all $a \in A$ and $x \in R$.  In this case, $1 \otimes 1$ is an identity in
$R \# A$.
\snl
In the general case, we have similar results (see e.g.\ the formulas in
proposition 5.9) but we have to consider the
multiplier algebra $M(R \# A)$ of $R \# A$.  This is possible because we have
seen that the product in $R \# A$ is non-degenerate.  First, we get the following
result.

\inspr{5.7} Proposition \rm
There exists a linear map $\pi : A \rightarrow M(R \# A)$
defined by
$$\align
\pi (a) (x^\prime \# a^\prime)  &=  \sum a_{(1)} x^\prime \# a_{(2)} a^\prime \\
(x^\prime \# a^\prime) \pi(a)  &= \, x \# a^\prime a
\endalign
$$
whenever $a,a^\prime \in A$ and $x \in R$.  This map is a unital algebra
homomorphism.

\snl\bf Proof : \rm
We first check that the above formulas give a well-defined element
$\pi(a)$ in $M(R \# A)$.  Indeed, for all $a, a^\prime, a^{\prime \prime} \in A$
and $x^\prime, x^{\prime \prime} \in R$ we get
$$\align
(x^{\prime \prime} \# a^{\prime \prime}) (\sum a_{(1)} x^\prime \# a_{(2)}
a^\prime) & = \sum x^{\prime \prime} (a^{\prime \prime}_{(1)}a_{(1)} x^\prime) \#
a^{\prime \prime}_{(2)} a_{(2)} a^\prime \\
& = (x \# a^{\prime \prime} a)(x^\prime \# a^\prime)
\endalign
$$
so that indeed we have the formula
$$(x^{\prime \prime} \# a^{\prime \prime}) (\pi (a) (x^\prime \# a^\prime)) =
((x^{\prime \prime} \# a^{\prime \prime}) \pi (a)) (x^\prime \# a^\prime).$$
That $\pi$ is a homomorphism follows easily from the formula $(x^\prime \#
a^\prime) \pi(a) = x^\prime \# a^\prime a$.
\snl
Let us now prove that $\pi$ is unital.  We clearly have $(R \# A) \pi(A) = R
\# A$ because $A^2 = A$.  On the other hand, because $\Delta (A) (1 \otimes A) =
A \otimes A$, we also get $\pi (A) (R \# A) = R \# A$.  This completes the
proof.
\hfill $\blacksquare$
\einspr

Similarly, we get the following.

\inspr{5.8} Proposition \rm
There exists a linear map $\pi : R \rightarrow M(R \# A)$
given by
$$\align
\pi (x) (x^\prime \# a^\prime) & = x x' \# a^\prime \\
(x^\prime \# a^\prime) \pi (x) & = \sum x'(a_{(1)} x) \# a_{(2)}
\endalign
$$
whenever $a \in A$ and $x,x' \in R$.  This map is an algebra homomorphism
of $R$ into $M(R \# A)$.  If $R^2 = R$, then $\pi$ is unital.

\snl\bf Proof : \rm
For all $a^\prime,a^{\prime \prime} \in A$ and $x,x^\prime,x^{\prime
\prime} \in R$ we get
$$\align
(x^{\prime \prime} \# a^{\prime \prime})(xx^\prime \# a^\prime)  & = \sum
x^{\prime \prime} (a^{\prime \prime}_{(1)} (xx^\prime)) \# a_{(2)} a^\prime \\
& = \sum x^{\prime \prime} (a^{\prime \prime}_{(1)} x) (a^{\prime \prime}_{(2)}
x^\prime) \# a_{(3)} a^\prime \\
& = (\sum x^{\prime \prime} (a^{\prime \prime}_{(1)} x) \# a^{\prime
\prime}_{(2)})(x^\prime \# a^\prime).
\endalign
$$
so that
$$(x^{\prime \prime} \# a^{\prime \prime})(\pi (x) (x^\prime \# a^\prime)) =
((x^{\prime \prime} \# a^{\prime \prime} ) \pi (x)) (x^\prime \# a^\prime)$$
and $\pi(x)$ is a well-defined element in $M(R \# A)$.
\snl
Again from the first formula it follows that $\pi$ is a homomorphism.
\snl
When $R^2 = R$ we get $\pi (R) (R \# A) = R \# A$.
\snl
On the other hand, if we replace $x^\prime$ by $a^\prime x^\prime$ and use that
$AR = R$ and that $\Delta (A) (A \otimes 1) = A \otimes A$, we get that $(R \#
A)\pi(R)$ is also equal to $R^2 \# A$.  So again, when $R^2 = R$ we get $(R \#
A) \pi (R) = R \# A$.
\hfill $\blacksquare$
\einspr

If we don't have $R^2 = R$, we just get $(R \# A) \pi (R) = \pi (R) (R \# A) =
R^2 \# A$. Now, in general we need unital homomorphisms into the
multiplier algebra in order to be able to extend them to multiplier
algebras. So, we have a slight problem here when $R^2\neq R$. We will
see after the next result that there is a way around it in this case.
\snl
We first prove some formulas that we had earlier in the case of algebras
with identity.

\inspr{5.9} Proposition \rm
When $\pi : A \rightarrow M(R \# A)$ and $\pi : R \rightarrow
M(R \# A)$ are defined as in the previous propositions, we have
$$\align
\pi (x) \pi(a) & = x \# a \\
\pi (a) \pi(x) & = \sum a_{(1)} x \# a_{(2)}
\endalign
$$
for all $a \in A$ and $x \in R$.

\snl\bf Proof : \rm
Take $a,a^\prime \in A$ and $x,x^\prime \in R$.  Then we get
$$\align
\pi (x) \pi (a) (x^\prime \# a^\prime) & = \pi(x) \sum a_{(1)} x^\prime \#
a_{(2)} a^\prime \\
& = \sum x (a_{(1)} x^\prime) \# a_{(2)} a^\prime \\
& = (x \# a)(x^\prime \# a^\prime)
\endalign
$$
and
$$\align
\pi(a) \pi(x) (x^\prime \# a^\prime) & = \pi(a) (xx^\prime \# a^\prime) \\
& = \sum a_{(1)} (xx^\prime) \# a_{(2)} a^\prime \\
& = \sum (a_{(1)} x)(a_{(2)} x^\prime) \# a_{(3)} a^\prime  \\
& = (\sum a_{(1)} x \# a_{(2)})(x^\prime \# a^\prime).
\endalign
$$
\hfill $\blacksquare$
\einspr

It follows that $R \# A = \pi (R) \pi (A) = \pi (A) \pi (R)$. This
implies e.g.\ that we can extend $\pi$ from $R$ to $M(R)$ (also when
$R^2\neq R$). Indeed, given
$m\in M(R)$ we define $\pi(m)$ in $M(R\# A)$ by
$\pi(m)(\pi(x)\pi(a))=\pi(mx)\pi(a)$ and
$(\pi(a)\pi(x))\pi(m)=\pi(a)\pi(xm)$. We will need this in section 7.
The map
$$x\ot a \to \sum a_{(1)}x \ot a_{(2)}$$
from $R\ot A$ to itself will play a fundamental r\^ole further. We see
here that it determines the commutation rules between elements of $\pi (A)$
and $\pi(R)$.
\snl
We now are ready to prove a universal property (related with these
commutation rules).

\iinspr{5.10} Proposition \rm
Let $(A,\Delta)$ be a regular multiplier Hopf algebra,
acting on an algebra $R$.  Suppose that $C$ is an algebra over $\Bbb C$ and that
there exists homomorphisms $\pi : A \rightarrow M(C)$ and $\pi : R \rightarrow
M(C)$ such that
$$\pi(a) \pi(x) = \sum \pi (a_{(1)} x)\pi(a_{(2)})$$
for all $a \in A$ and $x \in R$.  Then there is a homomorphism $\pi : R \# A
\rightarrow M(C)$ such that
$$\pi (x \# a) = \pi (x) \pi (a).$$

\snl\bf Proof : \rm
Define $\pi : R \# A \rightarrow M(C)$ by $\pi (x \# a) = \pi(x)
\pi(a)$.  To verify that $\pi$ is a homomorphism, take
$a,a^\prime \in A$ and $x,x^\prime \in R$.  We have
$$\align
\pi ((x \# a)(x^\prime \# a^\prime)) & =  \pi (\sum (x(a_{(1)} x^\prime) \#
a_{(2)}a^\prime)) \\
& = \sum \pi (x(a_{(1)} x^\prime)) \pi (a_{(2)} a^\prime) \\
& = \sum \pi (x) \pi (a_{(1)} x^\prime) \pi(a_{(2)}) \pi (a^\prime) \\
& = \pi(x) \pi(a)\pi(x^\prime)\pi(a^\prime) \\
& = \pi (x \# a)\pi(x^\prime \# a^\prime).
\endalign
$$
\hfill $\blacksquare$
\einspr

If we assume that $\pi$ is unital on $A$, i.e. $\pi (A)C = C$ and $C \pi
(A) = C$, then we can rewrite the condition.  Indeed, for all $a,a^\prime \in
A$ and $x \in R$ and  $b \in C$ we have
$$\align
\sum \pi(a_{(1)}) \pi (x) \pi (S(a_{(2)} ) a^\prime) b & = \sum \pi (a_{(1)} x)
\pi (a_{(2)} S(a_{(3)}) a^\prime) b \\
& = \pi (ax) \pi (a^\prime)b
\endalign
$$
and since $\pi(A)C = C$ we can interprete this formula as
$$\pi(ax) = \sum \pi (a_{(1)}) \pi (x) \pi (S(a_{(2)})).$$
The difficulty with the last expression is that non of the $a_{(k)}$ seem to be
covered.  This sum must be considered in the {\it strict topology} on
$M(C)$.
Indeed, if we multiply to the left or to the right with an element in
$B$, we obtain a covering for $a_{(1)}$ or for $a_{(2)}$ and we end up
with finite sums.
\snl
Now we look at a similar result for covariant modules.  Suppose we have
a left $R \# A$-module $V$ which is unital.  Then, we know that we can extend $V$ to a
left $M(R \# A)$-module.  Because we have our homomorphisms $\pi : A \rightarrow
M(R \# A)$ and $\pi : R \rightarrow M(R \# A)$, we have actions of $A$ and $R$
on $V$.  From the formula
$$\pi(a) \pi(x) = \sum \pi (a_{(1)} x) \pi (a_{(2)})$$
we get, in module formulation
$$a(xv) = \sum (a_{(1)} x) a_{(2)} v.$$
This precisely means that $V$ is a covariant $A$-$R$-module.  Remark that $V$ is a
unital $A$-module as $AV = \pi (A) (R \# A)V = (R \# A) V = V$.
\snl
We now prove the converse (which we announced already in the beginning
of this section).

\iinspr{5.11} Proposition \rm
Suppose that $V$ is a covariant $A$-$R$-module.  Then $V$ is
also a $R \# A$-module and the action of $R \# A$ is given by $(x \# a)v =
x(av)$ for all $x \in R$, $a \in A$ and $v \in V$.

\snl\bf Proof : \rm
Take $a,a^\prime \in A$, $x,x^\prime \in R$ and $v \in V$.  Then we get
$$\align
((x \# a)(x^\prime \# a^\prime))v & = (\sum x(a_{(1)} x^\prime) \# a_{(2)}
a^\prime)v \\
& = \sum (x(a_{(1)} x^\prime)) ((a_{(2)} a^\prime)v) \\
& = \sum x(a_{(1)} x^\prime) a_{(2)} a^\prime v \\
& = x(a(x^\prime a^\prime v)) \\
& = (x \# a)((x^\prime \# a^\prime)v)
\endalign
$$
and we get an action of $R \# A$ on $V$.
\hfill $\blacksquare$
\einspr

When is the action unital ?  If we assume that $V$ is unital, both for
$A$ and for $R$, then we get
$$(R \# A)V = R(AV) = RV = V.$$
Conversely, if $(R \# A)V = V$, then $RAV = V$ and $ARV = V$ and these imply
that $RV = V$ and $AV = V$.
\snl
Because $R$ is not assumed to have local units (in the sense of
proposition 2.2), we also want to look at the non-degeneracy of the
module actions (in the sense of proposition 3.2).
Suppose first that $V$ is a non-degenerate $R$-module. Let $v\in V$ be
given and assume that $(x\# a)v=0$ for all $a\in A$ and $x\in R$. This
means that $xav=0$ for all $x\in R$ and $a\in A$. Because $V$ is assumed
to be a non-degenerate $R$-module, it follows that $av=0$ for all $a\in
A$. Because $V$ is automatically a non-degenerate $A$-module, we get
$v=0$. Hence, $V$ is also a non-degenerate $(R\# A)$-module. Conversely,
assume that $V$ is a non-degenerate $(R\# A)$-module. Again take any
$v\in V$ and assume that $xv=0$ for all $x\in R$. Then $axv=0$ for all
$x\in R$ and $a\in A$. Using the twist map, we get also that $xav=0$ for
all $x\in R$ and $a\in A$. Hence, we must have $v=0$.
\snl
So, we get a one-to-one correspondence between unital left $R \# A$-modules
and covariant $R$-$A$-modules with unital actions. And if one is
non-degenerate, so is the other. This result is
similar to a well-known result for actions of locally compact groups on
C$^*$-algebras (see section 7.6 in [P]).
\abs
\head Inner and Equivalent Actions\endhead
\abs
We will first show that the smash product is trivial when the action is inner.
Recall that the action of $A$ on $R$ is called inner if there is a
unital homomorphism $\gamma : A \rightarrow M(R)$ such that
$$ax = \sum \gamma (a_{(1)})x \gamma (S(a_{(2)}))$$
for all $a \in A$ and $x \in R$ (cf. definition 4.13).

\iinspr{5.12} Proposition \rm
If the action of $A$ on $R$ is inner, then $R \# A$ is
isomorphic with $R \otimes A$, now considered with the usual tensor product
algebra structure.

\snl\bf Proof : \rm
Define $\varphi : R \# A \rightarrow R \otimes A$ by
$$\varphi (x \# a) = \sum x \gamma (a_{(1)}) \otimes a_{(2)}.$$
We first show that $\varphi$ is well-defined.  Because we assume $\gamma$ to be
unital, we can write $x$ as $\sum x_i \gamma (a_i)$.  Therefore, we see
that $a_{(1)}$ will be covered in the above definition.  Moreover, if $\sum x_i
\gamma (a_i) = 0$ then
$$\sum x_i \gamma (a_i a_{(1)}) \otimes a_{(2)} a^\prime = \sum x_i \gamma (a_i)
\gamma (a_{(1)}) \otimes a_{(2)} a^\prime = 0$$
for all $a^\prime$ so that also
$$\sum x_i \gamma (a_i a_{(1)}) \otimes a_{(2)} = 0.$$
This means that $\varphi$ is well-defined.
\snl
Similarly, we can define $\psi : R
\otimes A \rightarrow R \# A$ by
$$\psi (x \otimes a) = \sum x \gamma
(S(a_{(1)})) \otimes a_{(2)}.$$
\snl
An easy calculation will show that $\varphi \psi = \psi \varphi = \iota$.
\snl
Also, for all $a,a^\prime \in A$ and $x,x^\prime \in R$ we have
$$\align
\varphi ((x \# a)(x^\prime \# a^\prime)) & = \varphi(\sum x(a_{(1)} x^\prime) \#
a_{(2)} a^\prime) \\
& = \sum x (a_{(1)} x^\prime) \gamma (a_{(2)} a^\prime_{(1)}) \otimes a_{(3)}
a^\prime_{(2)} \\
& = \sum x \gamma (a_{(1)}) x^\prime \gamma (a^\prime_{(1)}) \otimes a_{(2)}
a^\prime_{(2)} \\
& = \varphi (x \# a) \varphi (x^\prime \# a^\prime).
\endalign
$$
This proves the result. \hfill $\blacksquare$
\einspr

Next, we will extend the above result to the case of equivalent actions.
We can only do this in the case of an action of a Hopf algebra because
we only have defined cocycle equivalent actions properly in that case. Recall
that two actions $\triangleright_1$ and $\triangleright_2$ are called (cocycle) equivalent (see
definition 4.14) if there is a linear map $\gamma : A \rightarrow M(R)$
satisfying $\gamma(1)=1$ and
\snl
\quad i) $\gamma(aa')=\sum \gamma (a_{(1})(a_{(2)}\triangleright_1
\gamma(a'))$
\snl
\quad ii) $\sum(a_{(1)}\triangleright_2 x)\gamma(a_{(2)})
=\sum\gamma(a_{(1)})(a_{(2)}\triangleright_1 x)$
\snl
for all $a,a'\in A$ and $x\in R$.
We get the following result.

\iinspr{5.13} Proposition \rm
If two actions $\triangleright_1$ and $\triangleright_2$ of a Hopf
algebra $A$ on $R$ are cocycle
equivalent, then the corresponding
smash products $R \#_1 A$ and $R \#_2 A$ are isomorphic.

\snl\bf Proof : \rm
Define $\varphi : R \#_2 A \rightarrow R \#_1 A$ by
$$\varphi (x \#_2 a) = \sum x \gamma (a_{(1)} ) \#_1 a_{(2)}$$
whenever $x \in R$ and $a \in A$.
Here, this map will be well-defined since we are working with a Hopf
algebra.  We
claim that $\varphi$ has an inverse $\psi$ given by
$$\psi (x \#_1 a) = \sum x(a_{(1)} \triangleright_1 \gamma(S(a_{(2)}))) \#_2
a_{(3)}.$$
\snl
We get
$$\align
\psi \varphi (x \#_2 a) & = \psi (\sum x \gamma (a_{(1)}) \#_1 a_{(2)}) \\
& = \sum x \gamma (a_{(1)}) (a_{(2)} \triangleright_1 \gamma
(S(a_{(3)}))) \#_2 a_{(4)}
\\
& = \sum x \gamma (a_{(1)} Sa_{(2)}) \#_2 a_{(3)} \\
& = x \#_2 a
\endalign
$$
because $\gamma(1)=1$. Similarly we obtain that
$$\align \varphi \psi (x \#_1 a)
  &=\varphi(\sum x(a_{(1)}\trr_1 \gamma(S(a_{(2)}))) \#_2 a_{(3)}) \\
  &=\sum x(a_{(1)}\trr_1 \gamma(S(a_{(2)}))\gamma(a_{(3)}))\#_1 a_{(4)}\\
  &=\sum
  x(a_{(1)}\trr_1(\gamma(S(a_{(3)}))(S(a_{(2)}))\trr_1\gamma(a_{(4)}))) \#_1
  a_{(5)}\\
  &=\sum x (a_{(1)}\trr_1\gamma(S(a_{(2)}a_{(3)}))) \#_1 a_{(4)} \\
  &=x\#_1 a.
  \endalign$$

\snl
Now, we verify that we have a homomorphism of algebras.  Take $x,x^\prime \in R$
and $a,a' \in A$.  Then
$$\align
\varphi ((x \#_2 a)(x^\prime \#_2 a^\prime)) & = \varphi (\sum x (a_{(1)}
\triangleright_2 x^\prime) \#_2 a_{(2)} a^\prime) \\
& = \sum x (a_{(1)} \triangleright_2 x^\prime) \gamma (a_{(2)} a^\prime_{(1)}) \#_1
a_{(3)} a^\prime_{(2)} \\
& = \sum x (a_{(1)} \triangleright_2 x^\prime) \gamma (a_{(2)})
(a_{(3)} \triangleright_1 \gamma (a^\prime_{(1)})) \#_1 a_{(4)} a^\prime_{(2)} \\
& = \sum x \gamma (a_{(1)} (a_{(2)} \triangleright_1 x^\prime)
(a_{(3)} \triangleright_1 \gamma
(a^\prime)) \#_1 a_{(4)} a^\prime_{(2)}.
\endalign
$$
On the other hand
$$\align
\varphi (x \#_2 a) \varphi (x^\prime \#_2 a^\prime) & = \sum (x \gamma (a_{(1)})
\#_1 a_{(2)})(x^\prime \gamma (a^\prime_{(1)}) \#_1 a^\prime_{(2)}) \\
& = \sum x \gamma (a_{(1)}) (a_{(2)} \triangleright_1 (x^\prime \gamma (a^\prime_{(1)} )))
\#_1 a_{(3)} a^\prime_{(2)}\\
& = \sum x \gamma (a_{(1)} (a_{(2)} \triangleright_1 x^\prime)
(a_{(3)} \triangleright_1 \gamma
(a^\prime)) \#_1 a_{(4)} a^\prime_{(2)}.
\endalign
$$
and this is precisely $\varphi ((x \#_2 a)(x^\prime \#_2 a^\prime))$.
\hfill$\blacksquare$
\einspr

As we have mentioned already in the previous section, where we
introduced this notion of equivalence, it will probably be not so easy
to extend this to the multiplier Hopf algebra case but it should be
possible, provided the right non-degeneracy condition on $\gamma$ is
assumed. In any case, the proof might become quite involved.
\abs
\abs
\head 6. Special Cases arising from Dual Pairs\endhead
\abs
In this section, we will consider some special cases, coming from
a dual pair of regular multiplier Hopf algebras (in the sense of
[D-VD]). In the case of the pairing between an algebraic quantum group
$A$ and its dual $\hat A$, we will be able to give an explicit (and
simple) form of the smash product. This case is close to the
finite-dimensional case.
\snl
The results that we obtain here will not only serve as an illustration
of the results of the previous sections, but will also be basic for the
duality, proved and discussed in the next section.
\snl
Let us first recall the notion of a pairing between two regular
multiplier Hopf algebas (as introduced in [D-VD]).
\snl
Start with two regular multiplier Hopf algebras $A$ and $B$ and a
non-degenerate bilinear from $\langle\ ,\ \rangle$ from $A\times B$ to
$\Bbb C$. For all $a\in A$ and $b\in B$ we can define an element
$a\triangleright b$ in $M(B)$ by
$$\align (a\trr b)b'&= \sum\langle a,b_{(2)}\rangle b_{(1)} b' \\
         b'(a\trr b)&= \sum\langle a,b_{(2)}\rangle b'b_{(1)}
         \endalign$$
whenever $b'\in B$. Let us now assume that the pairing is such that
$a\trr b\in B$ for all $a\in A$ and $b\in B$. Then, it makes sense to
require that
$$\langle a',a\trr b\rangle =\langle a'a,b\rangle$$
for all $a,a'\in A$ and $b\in B$. This essentially means that the
product in $A$ is dual to the coproduct in $B$. Then, for all
$a,a',a''\in A$ and $b\in B$ we will have
$$\align \langle a, (a'a'')\trr b\rangle
         &= \langle a(a'a''), b\rangle \\
         &= \langle aa', a''\trr b\rangle \\
         &= \langle a, a'\trr (a\trr b)\rangle \endalign$$
and because we have taken our pairing to be non-degenerate, we get
$$(a'a'')\trr b= a'\trr (a\trr b).$$
Therefore, $B$ is a left $A$-module.
\snl
This observation was the main motivation for the introduction of the
notion of a {\it prepairing} in [D-VD]. In fact, we impose four
conditions like the ones above. We require, for all $a\in A$ and $b\in
B$,
$$\align &\sum\langle a_{(1)}, b\rangle a_{(2)} \in A \\
         &\sum\langle a_{(2)}, b\rangle a_{(1)} \in A \\
         &\sum\langle a, b_{(1)}\rangle b_{(2)} \in B \\
         &\sum\langle a, b_{(2)}\rangle b_{(1)} \in B \endalign$$
and for all $a,a'\in A$ and $b,b'\in B$,
$$\align \langle \sum \langle a_{(1)},b \rangle a_{(2)},b'\rangle
            &=\langle a, bb'\rangle \\
         \langle \sum \langle a_{(2)}, b\rangle a_{(1)},b'\rangle
            &= \langle a, b'b\rangle \\
         \langle a',\sum \langle a, b_{(1)}\rangle b_{(2)}\rangle
            &= \langle aa', b\rangle \\
         \langle a',\sum \langle a, b_{(2)}\rangle b_{(1)}\rangle
            &=\langle a'a, b\rangle.\endalign$$
It follows that four modules are involved. We have that $B$ is a left
and a right $A$-module (in fact, an $A$-bimodule) and that $A$ is a left
and a right $B$-module (also a $B$-bimodule). We will denote the left
actions by $\trr$ and the right actions by $\trl$ (as in [D-VD]).
\snl
In proposition 2.8 of [D-VD], we have shown that all of these modules
are unital (in the sense of definition 3.1 of this paper) if only one of
them is unital. This is the case when the prepairing is a pairing
(definition 2.9 of [D-VD]).
\snl
Remark that the fact that these actions are unital also implies that
e.g.\ in the expression
$$\sum \langle a_{(1)},b\rangle a_{(2)},$$
the element $a_{(1)}$ is also {\it covered} by $b$ in the following
sense. Take an element $e\in A$ such that $e\trr b=b$ (cf.\ the remark
before proposition 3.3). Then we will have
$$\sum\langle a_{(1)},b\rangle a_{(2)}=\sum\langle a_{(1)}e,b\rangle
a_{(2)}$$
so that $a_{(1)}$ is covered in the usual sense. This covering of
$a_{(1)}$ by $b$ through the pairing in turn implies that
$$\sum \langle a_{(1)},b\rangle a_{(2)} \in A$$
for all $a\in A$ and $b\in B$. This shows that the axioms for a pairing
are not independent from each other. These observations are important
for the formulas and calculations further in this paper.
\nl
In what follows, we consider a dual pair $(A,B)$ of regular multiplier
Hopf algebras in the sense of definition 2.9 of [D-VD] (as discussed
before).
\snl
In the first proposition of this section, we show that $A$ is a
$B$-module algebra and that $B$ is a $A$-module algebra. The result
indicates once more that the notion of a pairing, as introduced in
[D-VD] is indeed a natural one (and not only for the treatment of the
quantum double).

\inspr{6.1} Proposition \rm
The left action of $A$ on $B$ defined by
$$a\trr b = \sum\langle a,b_{(2)}\rangle b_{(1)}$$
makes $B$ into a left $A$-module algebra (in the sense of definition
4.1).
\bf\snl Proof : \rm
We just have to show that
$$a\trr (bb') = \sum (a_{(1)}\trr b)(a_{(2)}\trr b')$$
for all $a\in A$ and $b,b'\in B$. Now, choose $e\in A$ so that $b=e\trr
b$. Take also $a'\in A$. Then
$$\align \langle a',\sum(a_{(1)}\trr b)(a_{(2)}\trr b')\rangle
     &=  \sum\langle a',((a_{(1)}e)\trr b)(a_{(2)}\trr b')\rangle  \\
     &= \sum_{(a)} \langle \sum_{(a')} \langle a'_{(1)},(a_{(1)}e)\trr b\rangle
        a'_{(2)},a_{(2)}\trr b'\rangle \\
     &= \sum \langle a'_{(1)}a_{(1)}e, b\rangle
        \langle a'_{(2)}a_{(2)}, b'\rangle \\
     &= \langle\sum  \langle(aa')_{(1)},b\rangle (aa')_{(2)}, b'\rangle \\
     &= \langle a'a, bb' \rangle = \langle a', a\trr (bb')\rangle.
\endalign$$
\hfill $\blacksquare$
\einspr

So, loosely speaking, $B$ is a left $A$-module because the product in
$A$ is dual to the coproduct on $B$ and $B$ is a left $A$-module algebra
because the product in $B$ is dual to the coproduct on $A$.
\snl
Because we have an action of $A$ on $B$, we can define the smash product
$B\# A$. Similarly, we have an action of $B$ on $A$ given by
$$b\trr a=\sum \langle a_{(2)},b \rangle a_{(1)}$$
whenever $a\in A$ and $b\in B$. This will yield another smash product
$A\# B$. We will now study these algebras and the relation between the
two.
\snl
First, observe that it is an immediate consequence of the definitions
that the product in $B\# A$ is given by the formula
$$(b\# a)(b'\# a')=\sum\langle a_{(1)},b'_{(2)}\rangle bb'_{(1)}\#
a_{(2)} a'.$$
Remark that there is no problem with this formula as $b'_{(1)}$ and
$a_{(2)}$ are both covered.
\snl
In the beginning of section 5,
we have seen that, in general, there is a natural action
of $R\# A$ on $R$ when $R$ is a left $A$-module algebra. The action is
given by $(x\# a)x'=x(ax')$ whenever $a\in A$ and $x,x'\in R$. This
action need not be faithful. In this special case however, we do have a
faithful action:

\inspr{6.2} Proposition \rm
The space $B$ is made into a left $(B\# A)$-module by
$$(b\# a)b'=b(a\trr b')$$
whenever $a\in A$ and $b,b'\in B$. This action of $B\# A$ is faithful.
\snl\bf Proof : \rm
It is straightforward to check that $B$ is a left $(B\# A)$-module for
the action defined above. This was shown already in section 5.
To show that the
module is faithful, let $\sum b_i\# a_i \in B\# A$ and assume that
$$\sum (b_i\# a_i)b'=\sum\langle a_i,b'_{(2)}\rangle b_ib'_{(1)} = 0$$
for all $b'\in B$. If we multiply with $b''$ on the right and if we use
that $\Delta(B)(B\ot 1)=B\ot B$, we get that
$$\sum\langle a_i,b'\rangle b_ib''=0$$
for all $b',b''\in B$. We can cancel $b''$ because the product in $B$ is
non-degenerate and using the non-degeneracy of the pairing, we can
conclude that $\sum b_i\ot a_i=0$ in $B\ot A$. This completes the proof.
\hfill $\blacksquare$
\einspr

Compare this result with lemma 9.4.2 of [M2] where the faithfulness of
the map $\lambda$ follows from the bijectivity of the antipode.
\snl
Related with this proposition is the following characterization of the
algebra $B\# A$.

\inspr{6.3} Proposition \rm
Consider the algebra $C$ generated by the algebras $A$ and $B$ subject
to the following commutation relations
$$ab=\sum\langle a_{(1)},b_{(2)}\rangle b_{(1)}a_{(2)}$$
whenever $a\in A$ and $b\in B$. Then $b\# a\to ba$ gives an injective
homomorphism of $B\# A$ into $C$.
\snl\bf Proof : \rm
First remark that
$$\sum \langle a_{(1)}, b_{(2)}\rangle a_{(2)} \ot b_{(1)}$$
is defined within $A\ot B$ (see proposition 2.8 of [D-VD]). It follows
that the algebra $C$ is well-defined. We have an action of $A$ and an
action of $B$ on $B$. Moreover, we have
$$\align
   a\trr (bb')&=\sum(a_{(1)}\trr b)(a_{(2)}\trr b')\\
              &=\sum\langle a_{(1)},b_{(2)}\rangle b_{(1)} (a_{(2)}\trr b')
\endalign$$
for all $a\in A$ and $b,b'\in B$. But this precisely means that these
two actions of $A$ and $B$ satisfy the defining commutation rules for
$C$ and therefore, these two actions yield an action of $C$ on $B$.
Now, since $B$ is a faithfull $(B\# A)$-module for the action given by
$(b\# a)b'=b(a\trr b')$, we see that the map $b\# a \to ba$ gives indeed
an injective homomorphism of $B\# A$ into $C$.
\hfill $\blacksquare$ \einspr

Remark that the result above is naturally related with the universal
property of $R\# A$ as proved in 5.10.
\snl
The commutation rules in 6.3 can also be expressed as
$$ba=\sum\langle S^{-1} a_{(1)},b_{(2)}\rangle a_{(2)}b_{(1)}.$$
This can easily be verified and in fact, it follows from the fact that
the maps
$$\align a\ot b &\to \sum \langle a_{(1)},b_{(2)}\rangle a_{(2)}\ot b_{(1)}\\
         a\ot b &\to \sum \langle S^{-1}a_{(1)},b_{(2)}\rangle a_{(2)}\ot b_{(1)}
\endalign$$
are each others inverses (see also proposition 2.8 in [D-VD]).
\snl
So, $B\# A$ can be considered as the span of the elements $ab$ in the
algebra generated by $A$ and $B$ subject to the commutation rules above.
Similarly, the algebra $A\# B$ can be characterized as the span of the
elemens $ab$ in the algebra generated by $A$ and $B$ subject to the
commutation relations
$$ba=\sum\langle a_{(2)},b_{(1)}\rangle a_{(1)}b_{(2)}$$
whenever $a\in A$ and $b\in B$.
\snl
The commutation rules above are of the same nature as the Heisenberg
commutation rules (e.g.\ in the case of the pairing between an abelian
group $G$ and its dual group $\hat G$). Therefore, the algebra $C$ in
this proposition is sometimes called the Heisenberg double (cf.\ [M2]).
The quantum double  (of Drinfel'd) on the other hand can be
characterized in a similar way, but by using different commutation rules
(see [D-VD]). In the first case, the commutation rules are related with
the Pentagon equation while in the second case, they are related with
the Yang-Baxter equation (see e.g.\ [VD2]).
\snl
Using proposition 6.3,
it becomes very easy to get the following relation between $B\# A$
and $A\# B$.

\inspr{6.4} Proposition \rm
The map $b\# a\to S^{-1}a\# Sb$ defines a anti-isomorphism of $B\#
A$ with $A\# B$.
\snl\bf Proof : \rm
It is straighforward to check this result directly. On the other hand,
one can also easily verify that this map will convert the commutation
rules defining $B\# A$ into those defining $A\# B$. Indeed, for all
$a\in A$ and $b\in B$ we get
$$\sum\langle(S^{-1}a)_{(2)}, (Sb)_{(1)}\rangle (S^{-1}a)_{(1)}\ot
(Sb)_{(2)}
=\sum\langle a_{(1)}, b_{(2)}\rangle S^{-1}a_{(2)}\ot Sb_{(1)}.$$
\hfill $\blacksquare$ \einspr

Compare with proposition 2.14 of [D-VD].
\snl
We can use this anti-isomorphism of $B\# A$ with $A\# B$ to obtain a
standard right action of $B\# A$ on $A$. Indeed, we can make $A$ into a
right $(B\# A)$-module by letting
$$\align a'(b\# a)&=(S^{-1}a\# Sb)a'\\
                  &=(S^{-1}a)((Sb)\trr a')\\
                  &=\sum\langle a'_{(2)},Sb\rangle (S^{-1}a)a'_{(1)}.
\endalign$$
If we replace $a'$ by $S^{-1}a'$, we find
$$(S^{-1}a')(b\# a)=\sum\langle a'_{(1)},b \rangle S^{-1}(a'_{(2)}a)$$
and so we get a modified $(B\# A)$-module structure on $A$ given by
$$a'(b\# a)=(a'\trl b)a$$
whenever $a,a'\in A$ and $b\in B$. Of course also this action is
faithful. This corresponds with the faithfulness of $\rho$ in lemma
9.4.2 of [M2].
\snl
In the next section, where we discuss the duality for actions, we will
come back to these two $B\# A$ module structures.
\snl
We will specialize to the pairing of an algebraic quantum group $A$
with its dual $\hat A$. But before we do this, we will prove a result
about fixed points that will be used in the next section.
\snl
\inspr{6.5} Proposition \rm
The only fixed points in $M(B)$ for the action of $A$ are multiples of
the identity.
\snl\bf Proof : \rm
Let $m$ be in $M(B)$ and a fixed point for the action of $A$. We have
seen in proposition 4.11 that then $$a\trr (mb)=m(a\trr b)$$ for all $a\in
A$ and $b\in B$. If we apply the counit to this equation, we get
$$\langle a, mb \rangle = \sum \epsilon (mb_{(1)}) \langle a, b_{(2)}
\rangle$$
for all $a$ and $b$. Then we get
$$ mb  = \sum \epsilon (mb_{(1)}) b_{(2)}$$
for all $b$.
Therefore we have for all $b,b'\in B$ that
$$\align \epsilon(mb)b'&=\sum\epsilon(m b_{(1)}) b_{(2)} S(b_{(3)}) b'\\
          &=\sum (mb_{(1)}) (S(b_{(2)}) b' \\
          &=\sum mb_{(1)}S(b_{(2)})b'\\
          &= \epsilon(b) mb'.
  \endalign$$
So, we get $m=\epsilon(mb)1$ when $b$ is chosen so that $\epsilon(b)=1$.
\hfill $\blacksquare$
\einspr

Remark that in the proof above, we only have used that $m$ is a left
multiplier. In fact, it is for such multipliers that we will need this result in
the next section.
\nl
Now, in what follows, $A$ is an algebraic quantum group (cf.\ definition
2.4) and $\hat A$ denotes the dual of $A$ (cf.\ theorem 2.5). The
duality between $A$ and $\hat A$ makes $(A,\hat A)$ into a dual pair
satisfying the axioms mentioned in the beginning of this section (see
[D-VD]). So, we
can apply the previous results with $B$ replaced by $\hat A$. In
particular we have the left action of $\hat A$ on $A$ and the
corresponding smash product $A\# \hat A$. We also have the action of
$A\# \hat A$ on $A$ given by
$$(a\# b)a'=\sum\langle a'_{(2)},b \rangle aa'_{(1)}$$
where now $a,a'\in A$ and $b\in \hat A$.
\snl
We get the following result.

\inspr{6.6} Proposition \rm
The algebra $A\# \hat A$, as acting on $A$, is the algebra of linear
combinations of the (rank one) maps from $A$ to $A$ of the form $a'\to
\langle a',b\rangle a$ where $a\in A$ and $b\in \hat A$.
\snl\bf Proof : \rm
Let $c\in A$ and consider $b=\varphi(c\,\cdot\,)$ in $\hat A$ where
$\varphi$ is a given left integral on $A$. Then
$$\align \sum\varphi(c a'_{(2)}) aa'_{(1)}
         &=\sum\varphi(c_{(3)} a'_{(2)}) a S(c_{(1)})c_{(2)}a'_{(1)}\\
         &=\sum\varphi(c_{(2)} a') a S(c_{(1)}).
\endalign$$
It follows from the fact that $\Delta(A)(A\ot 1)=A\ot A$ and from the
bijectivity of the antipode that the algebra $A\# \hat A$, as acting on
$A$, is precisely the linear span of the operators on $A$ of the form
$a'\to \varphi(ca')a$ with $a,c\in A$.
\hfill $\blacksquare$
\einspr

It is not hard to obtain the structure of $A\# \hat A$ from this.
Consider $A\ot \hat A$, together with the product
$$(a\ot b)(a'\ot b')=\langle a',b\rangle a\ot b'$$
whenever $a,a'\in A$ and $b,b'\in \hat A$. When $A\ot\hat A$ is considered
with this product, we will use $A\diamondsuit \hat A$ and for the
elements $a\diamondsuit b$. This algebra has an obvious faithful action
on $A$ given by
$$(a\diamondsuit b)a'=\langle a',b\rangle a$$
whenever $a,a'\in A$ and $b\in \hat A$. Then the following result is an
immediate consequence of proposition 6.5.

\inspr{6.7} Proposition \rm
The map $\gamma : A\# \hat A \to A \diamondsuit \hat A$, defined by
$$\gamma (a\# \varphi(c\,\cdot\,))=
    \sum a S(c_{(1)}) \diamondsuit \varphi(c_{(2)}\,\cdot\,)$$
is an isomorphism of algebras.
\einspr

Remark that the algebra structure of $A\diamondsuit \hat A$ no longer
depends on the products, nor on the coproducts of $A$ and $\hat A$. It
only depends on the bilinear form between the two vector spaces $A$ and
$\hat A$.
\snl
In a similar way, we can consider the algebra $\hat A\diamondsuit A$. We
get an obvious algebra isomorphism of $\hat A \# A$ with $\hat A
\diamondsuit A$.
\snl
We also have seen in proposition 6.4
that $A\# \hat A$ is anti-isomorphic with $\hat A \#
A$. It is not hard to see that also $A\diamondsuit \hat A$ and $\hat A
\diamondsuit A$ are anti-isomorphic. The flip will realize such an
anti-isomorphism.
\snl
In fact, it is also possible to construct an isomorphism of
$A\diamondsuit \hat A$ to $\hat A\diamondsuit A$. Any bijective linear
map from $A$ to $\hat A$ will do. It is also not so difficult to realize
an isomorphism from $A\# \hat A$ to $\hat A\# A$. Indeed, we know that
both algebras act on $A$. We have the natural left action of $A\# \hat
A$ considered before and given by
$$(a\# b)a'=\sum\langle a'_{(2)} ,b\rangle aa'_{(1)}$$
but we also have the right action of $\hat A\# A$ on $A$ given by
$$a'(b\# a)=\sum \langle a'_{(1)}, b \rangle a'_{(2)} a.$$
If now, we take $b=\psi(\,\cdot\, c)$ with $c\in A$ and a right integral
$\psi$ on $A$, we will find
$$\align \sum\psi(a'_{(1)} c) a'_{(2)} a
     &=\sum\psi(a'_{(1)}c_{(1)})a'_{(2)}c_{(2)}S(c_{(3)})a\\
     &=\sum\psi(a'c_{(1)})S(c_{(2)})a.
\endalign$$
As in the proof of 6.6, it follows again that the right action of $\hat
A\# A$ on $A$ is spanned by the rank one operators from $A$ to $A$ of
the form $a'\to \psi(a'c)a$ with $a\in A$ and $c\in A$. But also
$\hat A$ is the set of functionals on $A$ given by $\psi(\,\cdot\, c)$
with $c\in A$ and hence, we get again the algebra $A\diamondsuit \hat
A$.
\snl
Let us finish this section with a few remarks.
\snl
We have started with a dual pair $(A,B)$ of multiplier Hopf algebras and
we have assumed that the pairing is non-degenerate. In the case of a
degenerate pairing, similar results can be obtained, but one has to be
somewhat more careful.
\snl
A particular case is obtained when the algebra $A$ is actually a Hopf
algebra. Then, for $B$ one can take a Hopf subalgebra of the dual
algebra $A'$ of $A$. This is in fact the setting for the duality
considered by Montgomery in [M2]. If $B$ is large enough (i.e.\
separates points), then we have a non-degenerate pairing. Otherwise,
this is not the case.
\snl
Finally, the finite dimensional case is obviously a special case of the
pairing between $A$ and $\hat A$ for an algebraic quantum group. Then,
the smash product $\hat A\# A$ is simply the algebra of $n\times n$
matrices over $\Bbb C$ where $n$ is the dimension of $A$.
\FF
\head 7. Duality\endhead
\abs
In this section, we will obtain the main result of this paper. We first
define the dual action on the smash product and then we will study the
bismash product. We do this for any dual pair of regular multiplier Hopf
algebras. However, the nicest duality theorem is obtained for an
algebraic quantum group $A$. In this case, the bismash product is an
algebra with a very simple structure, no longer depending on the product
and the coproduct of $A$, nor on the action of $A$. This situation is
similar to (and in fact strongly connected with) the result obtained in
proposition 6.6 on the structure of $A\# \hat A$. The result
that we can prove in the case of a general pairing requires some extra
conditions which makes it really a result on coactions. Since we are
preparing a separate paper on coactions, we will only state the result
here and make some comments. In [VD-Z3], we will treat this result in
full detail. It should be noted however that the extra conditions needed
in the the case of a general pairing are automatically fulfilled for the
pair $(A,\hat A)$. So the duality theorem that we will prove in [VD-Z3]
is an extension of the result that we have here on $(A,\hat A)$.
\snl
So, as before,
let $(A,\Delta)$ be a regular multiplier Hopf algebra, acting on an
algebra $R$. Consider the smash product $R\# A$ as
defined and studied in section 5. Let $(B,\Delta)$ be another regular
multiplier Hopf algebra, paired with $(A,\Delta)$ in the sense of
[D-VD], as explained in the previous section. Again assume that the
pairing is non-degenerate.
\snl
Recall that we have a left action of $B$ on $A$ given by
$$b\trr a=\sum \langle a_{(2)}, b \rangle  a_{(1)}$$
and remark that in this formula, $b$ covers $a_{(2)}$ through the
pairing. This left action of $B$ on $A$ gives the {\it dual action} of
$B$ on $R\# A$ in the following definition.

\inspr{7.1} Definition \rm
Define a linear map from $B\ot (R\# A)$ to $R\# A$ by
$$b(x\# a)=x\#(b\trr a)=\sum\langle a_{(2)}, b\rangle\, x\# a_{(1)}$$
whenever $a\in A$, $b\in B$ and $x\in R$.
\einspr

It easily follows that this map makes $R\# A$ into a left $B$-module. It
is a unital module because $A$ is a unital $B$-module. But just as $A$
is a $B$-module algebra, the same is true for $R\# A$:

\inspr{7.2} Proposition \rm
The action defined above makes $R\# A$ into a $B$-module algebra.
\snl\bf Proof : \rm
Take $a,a'\in A$, $b\in B$ and $x,x'\in R$. Then
$$\align b((x'\# a')(x\# a))
    &= \sum b(x'(a'_{(1)}x)\# a'_{(2)}a)\\
    &= \sum x'(a'_{(1)}x) \# (b\trr (a'_{(2)}a))\\
    &= \sum x'(a'_{(1)}x) \# (b_{(1)}\trr a'_{(2)})(b_{(2)}\trr a)\\
    &= \sum (x'\# (b_{(1)}\trr a'))(x\# (b_{(2)}\trr a))\\
    &= \sum (b_{(1)}(x'\# a'))(b_{(2)}(x\# a)).\endalign$$
This proves that we have a $B$-module algebra.
\hfill $\blacksquare$
\einspr

Remark that we always have the appropriate coverings in the above
calculations. We e.g.\ have that $a'_{(1)}$ is covered by $x$ and that
$b_{(2)}$ is covered by $a$ or by $x\# a$. We also have used that the
left action of $B$ on $A$ commutes with the comultiplication in the
sense that
$$\Delta(b\trr a)(a'\ot 1)=\sum a_{(1)}a' \ot (b\trr a_{(2)})$$
for all $a,a'\in A$ and $b\in B$.
\snl
It is interesting to look at the fixed points in $M(R\# A)$ for the dual
action. Remember that we have a homomorphism $\pi : R \to M(R\# A)$
given by
$$\pi(x)(x'\# a')=xx'\# a'$$
(see proposition 5.8). It is not so difficult to see that $\pi(x)$ is a
fixed point. Indeed, if $x,x'\in R$, $a'\in A$ and $b\in B$, we get
$$b(\pi(x)(x'\# a'))=xx'\# b\trr a'=\pi(x)b(x'\# a')$$
and this means that $\pi(x)$ is a fixed point (see the remark following
proposition 4.11). As we have
mentioned already in section 5 (cf.\ the remark after proposition 5.9),
we can extend the homomorphism $\pi$ to the
multiplier algebra $M(R)$ and the formula is still given by
$$\pi(m)(x\# a)=mx\# a$$
whenever $a\in A$, $x\in R$ and $m\in M(R)$. It is again easy to see that
also $\pi(m)$ is a fixed point for the dual action.
\snl
In fact, we can prove the converse:

\inspr{7.3} Proposition \rm
The algebra of fixed points in $M(R\# A)$ for the dual action is
precisely the image of $M(R)$ under the canonical imbedding $\pi$.

\snl\bf Proof : \rm
So, let $z$ be a multiplier of $R\# A$ and assume that $z$ is invariant
under the dual action. Fix $x\in R$ and let $\omega\in R'$. Define a
linear map $\ell : A \to A$ by
$$\ell (a)=(\omega \ot \iota)(z(x\# a)).$$
For all $a,a'\in A$ we have
$$\ell(aa')=(\omega\ot\iota)(z(x\# aa'))
           =(\omega\ot\iota)(z(x\# a)\pi(a'))$$
and since $(x''\# a'')\pi(a')=x''\# a''a'$, we get
$$(\omega\ot\iota)(z(x\# a)\pi(a'))=((\omega\ot\iota)(z(x\# a)))a'$$
so that $\ell(aa')=\ell(a)a'$. This means that $\ell$ is a left multiplier of
$A$. We will write $\ell a$ for $\ell(a)$ in what follows.
\snl
Now, we will use that $z$ is fixed under the dual action. This means
that $b(z(x\# a))=z(b(x\# a))=z(x\# b\trr a)$. Because $b$ acts on the
second factor, we get
$$\align
b\trr \ell a&=(\omega\ot\iota)(b(z(x\# a)))    \\
          &=(\omega\ot\iota)(z(x\# b\trr a)) \\
          &=\ell (b\trr a).\endalign$$
So, we see that $\ell$ is a left multiplier on $A$ and fixed under the action
of $B$. We have shown in section 6 (remark after proposition 6.5)
that this implies that $\ell$ is a scalar
multiple of the identity. It follows that there is a number $\lambda$
such that
$$(\omega\ot\iota)(z(x\# a))=\lambda a$$
for all $a\in A$. This implies that for all $a\in A$ and $x\in R$ there
exists an element $y\in R$ such that
$$z(x\# a) = y\# a.$$
We will now show that $y$ does not depend on $a$. Take $a,a'\in A$.
Choose $e\in A$ so that $ea=a$ and $ea'=a'$. Take $y$ so that
$z(x\# e)=y\# e$. If we multiply this equation with $\pi(a)$ on the
right, we find that also $z(x\# a)=y\# a$. Similarly, we get $z(x\#
a')=y\# a'$. This implies that $y$ does not depend on $a$.
\snl
So, we get a linear map, again denoted by $\ell$ from $R$ to $R$ so that
$z(x\# a)=\ell(x)\# a$ for all $a\in A$ and $x\in R$. Using once more that
$z$ is a left multiplier, we obtain
$$\align \sum \ell(x) a_{(1)}x'\# a_{(2)}a'
       &=(\ell(x)\# a)(x'\# a') \\
       &=(z(x\# a))(x'\# a')    \\
       &=z((x\# a)(x'\# a'))    \\
       &=\sum\ell(x(a_{(1)}x'))\# a_{(2)}a'\endalign$$
for all $a,a'\in A$ and $x,x'\in R$.
If we apply $\epsilon$ to the second factor, we get
$$\ell(x(ax'))=\ell(x)ax'$$
for all $a\in A$ and $x,x'\in R$. Hence we have also
$\ell(xx')=\ell(x)x'$ so that $\ell$ is a left multiplier of $R$. We
will again in what follows denote $\ell(x)$ by $\ell x$.
\snl
Now, we will use the fact that $z$ is also a right multiplier to prove
that also $\ell$ is a right multiplier. For all $a\in A$, $x\in R$ and
$y\in R\# A$ we get
$$(yz)(x\# a)=y(\ell x\# a).$$
Take for $y$ an element of the form $\pi(a')\pi(x')$ with $a'\in A$ and
$x'\in R$ and write
$$\pi(a')\pi(x')z=\sum\pi(a_i)\pi(x_i)$$
to get
$$\pi(a')\pi(x'(\ell x))\pi(a)=\sum\pi(a_i)\pi(x_ix)\pi(a).$$
We can cancel $\pi(a)$ and use the injectivity of the map
$x\ot a\to \pi(a)\pi(x)$
on $R\ot A$ (because $\pi(a)\pi(x)=\sum a_{(1)}x\# a_{(2)}$) to obtain
$$x'(\ell x)\ot a'=\sum x_ix\ot a_i.$$
Now apply any linear functional on the second factor. Then we get for
every element $x'\in R$ and element $y$ in $R$ satisfying $x'(\ell x)=
yx$ for all $x\in R$. This implies that $\ell$ is in fact a two-sided
multiplier. This completes the proof.   \hfill $\blacksquare$
\einspr

Similar results are known for actions of locally compact quantum groups
on operator algebras (see e.g.\ proposition 4.12 in [VD1] for one of
them).
\nl
Now, we consider the {\it bismash product} $(R\# A)\# B$. The product of
two elements $(x\# a)\# b$ and $(x'\# a')\# b'$ in $(R\# A)\# B$ is
given by
$$\align ((x\# a)\# b)((x'\# a')\# b')
   &= \sum ((x\# a)(x'\# (b_{(1)}\trr a'))) \# b_{(2)}b' \\
   &= \sum \langle a'_{(2)},b_{(1)}\rangle ((x\# a)(x'\# a'_{(1)}))\#
      b_{(2)}b' \\
   &= \sum\langle  a'_{(2)},b_{(1)} \rangle ((x(a_{(1)}x'))\#
      a_{(2)}a'_{(1)})\# b_{(2)} b'.
  \endalign$$
Remark that in this last expression, $b'$ covers $b_{(2)}$ and $x'$
covers $a_{(1)}$. This implies that indirectly, also both $a'_{(1)}$ and
$a'_{(2)}$ are covered.
\snl
We have seen that, in general, $R$ can be made into a left $(R\#
A)$-module (cf.\ section 5). Here, we can make $R\# A$ into a left $((R\# A)\#
B)$-module. We have the following formula. Whenever $(x\# a)\# b$ is an
element in $(R\# A)\# B$ and $x'\# a'$ is in $R\# A$, we get
$$\align ((x\# a)\# b)(x'\# a')
   &= (x\# a)(x'\# (b\trr a'))\\
   &= \sum \langle a'_{(2)}, b \rangle (x\# a)(x'\# a'_{(1)})\\
   &= \sum \langle a'_{(2)}, b \rangle x(a_{(1)}x') \# a_{(2)}a'_{(1)}.
   \endalign$$
Remark again that $a_{(1)}$ is covered by $x'$ and $a'_{(2)}$ is
covered by $b$.
\snl
Just as in the case of the action of $A\# B$ on $A$ (see proposition 6.2),
also here
this action turns out to be faithful:

\inspr{7.4} Proposition \rm
$R\# A$ is a faithful $((R\# A)\# B)$-module for the action defined
above.
\snl\bf Proof : \rm
Assume that $\sum (x_i\# a_i)\# b_i$ is in $(R\# A)\# B$ and that
$$(\sum (x_i\# a_i) \# b_i)(x'\# a')=0$$
for all $x'\in R$ and $a'\in A$.
Then, using the above formula, we get
$$\sum\langle a'_{(2)}, b_i\rangle x_i(a_{i(1)}x')\# a_{i(2)}a'_{(1)} =0$$
for all $x'\in R$ and $a'\in A$. As in the proof of 6.2, multiply the
second factor in this smash product with $a''$ from the right and use
that $\Delta(A)(A\ot 1)=A\ot A$ to get
$$\sum\langle a', b_i\rangle x_i(a_{i(1)}x')\# a_{i(2)}a'' =0$$
for all $x'\in R$ and $a',a''\in A$. Then we can use the
non-degeneracy of the product in $R\# A$ to obtain that
$$\sum\langle a', b_i\rangle x_i \# a_i =0$$
for all $a'\in A$. Finally, because we assume the pairing to be
non-degenerate, we get $\sum (x_i\# a_i)\# b_i =0$.
\hfill $\blacksquare$
\einspr

Now, we will construct an equivalent module. We will denote the
underlying space again by $R\ot A$. Define a linear map $W: R\ot A \to
R\# A$ by
$$W(x\ot a)=\sum a_{(1)}x\# a_{(2)}.$$
This is well-defined because $a_{(1)}$ is covered by $x$. Remark that
$W$ is surjective because $\Delta(A)(A\ot 1)=A\ot A$ and $R$ is unital.
In fact, $W$ is a bijection and the inverse is given by
$$W^{-1}(x\# a)=\sum S^{-1}(a_{(1)})x\ot a_{(2)}$$
whenever $a\in A$ and $x\in R$. We have used a similar map in section 5.
In fact, it is the twisting map defining the commutation rules between
$\pi(A)$ and $\pi(R)$ (cf.\ proposition 5.10).
\snl
Now, we can use this vector space isomorphism $W$ to construct a new
action of $(R\# A)\# B$ on $R\ot A$. It will still be faithful and we
will use it to examine the bismash product $(R\# A)\# B$.
\snl
In the next proposition, we obtain the formulas giving this action.

\inspr{7.5} Proposition \rm
For $a,a'\in A$, $b\in B$ and $x,x'\in R$, we get
$$\align W^{-1}(x\# a) W (x'\ot a')
     &= \sum ((S^{-1}a'_{(1)})(S^{-1}a_{(1)})x)x'\ot a_{(2)} a'_{(2)}\\
     W^{-1} b W (x'\ot a')
     &= x'\ot (b\trr a').
   \endalign$$
\snl\bf Proof : \rm
We have
$$\align
     W^{-1}(x\# a) W (x'\ot a')
       &= W^{-1} \sum (x\# a)(a'_{(1)}x' \# a'_{(2)}) \\
       &= W^{-1} \sum x(a_{(1)}a'_{(1)}x') \# a_{(2)}a'_{(2)} \\
       &= \sum S^{-1} (a_{(2)}a'_{(2)})(x(a_{(1)}a'_{(1)}x'))
          \ot a_{(3)}a'_{(3)} \\
       &= \sum (S^{-1} (a_{(1)}a'_{(1)})x)x'\ot a_{(2)}a'_{(2)}
   \endalign$$
and
$$\align
     W^{-1} b W (x'\ot a')
        &= W^{-1} \sum a'_{(1)}x' \# (b\trr a'_{(2)}) \\
        &= W^{-1} \sum \langle  a'_{(3)}, b\rangle a'_{(1)}x' \# a'_{(2)}\\
        &= \sum \langle  a'_{(2)},b\rangle x'\ot a'_{(1)}.
   \endalign$$
\hfill $\blacksquare$
\einspr

If we combine these two formulas, we get
$$W^{-1} ((x\# a)\# b) W(x'\ot a')
     = \sum \langle a'_{(3)}, b \rangle ((
     S^{-1}a'_{(1)})(S^{-1}a_{(1)})x)x' \ot a_{(2)}a'_{(2)}$$
for all $a,a'\in A$, $b\in B$ and $x,x'\in R$.
\snl
We are interested in the algebra of linear maps on $R\ot A$ spanned by
such operators. For convenience, we will call this algebra $P$ in what
follows.
\snl
Now, we have seen that the map
$$x\ot a \to \sum S^{-1}(a_{(1)}) x \ot a_{(2)}$$
is a bijection of $R\ot A$ (in fact, it is the map $W^{-1}$).
It follows that the algebra $P$ is the span of all operators of the form
$$x'\ot a' \to \sum\langle a'_{(3)}, b \rangle ((S^{-1}a'_{(1)})x)x'
 \ot a a'_{(2)}$$
where still $a\in A$, $b\in B$ and $x\in R$.
We can recognise the action of $A\# B$ on $A$ in the last factor.
Indeed, the above operator can be rewritten as
$$x'\ot a' \to \sum ((S^{-1}a'_{(1)})x)x' \ot (a\# b) a'_{(2)}.$$
\snl
We will come back to this expression later.
\snl
Now, we first want to
specialize to the case of an algebraic quantum group.
\snl
So, in what follows, let $A$ be an algebraic quantum group and let $B$
be the dual $\hat A$ of $A$. Of course, we consider the pairing given by
the duality. We have seen in the previous section that the algebra $A\#
\hat A$, as acting on $A$ is the span of rank one operators of the form
$a' \to \langle a', b \rangle a$ where $a\in A$ and $b \in \hat A$ (see
proposition 6.6). Remember that we have used the notation $A \diamondsuit
\hat A$ for this algebra of finite rank operators.
\snl
We will now use the same technique as in the proof of proposition 6.6 to
obtain the following duality theorem.

\inspr{7.6} Theorem \rm
If $A$ is an algebraic quantum group, acting on an algebra $R$ and if
$\hat A$ is the dual of $A$, acting on the smash product $R\# A$ by
means of the dual action, then the bismash product $(R\# A)\# \hat A$ is
isomorphic with $R\ot (A\diamondsuit \hat A)$.
\snl \bf Proof : \rm
Let $a\in A$ and $b\in \hat A$. Write $b$ in the form
$\varphi(c\,\cdot\,)$ where $\varphi$ is a left integral on $A$ and
$c\in A$. As we had in the proof of proposition 6.6, we get
$$\sum \langle a'_{(2)}, b \rangle a'_{(1)}=
   \sum \varphi (c_{(2)} a') S(c_{(1)}).$$
So, we find that the algebra $P$ is spanned by operators of the form
$$x'\ot a'\to \sum \varphi(c_{(3)}a')(c_{(2)}x)x' \ot a S(c_{(1)})$$
where $a,c\in A$ and $x\in R$.
Because $(A\ot 1)\Delta(A)=A\ot A$, we get all operators of the
form
$$x'\ot a' \to \sum \varphi(c_{(2)}a')(c_{(1)}x)x'\ot a.$$
and by the fact that
$$x\ot c \to \sum c_{(1)}x\ot c_{(2)}$$
is surjective on $R\ot A$, we have that $P$ is spanned by operators of
the form
$$x'\ot a' \to xx' \ot \varphi (ca')a,$$
still with $x\in R$ and $a,c\in A$. Since we are dealing with faithful
actions of these two algebras, we get an isomorphism.
This completes the proof.
\hfill $\blacksquare$
\einspr

If we examine the proof more carefully and if we follow the different
steps in it, we can find an explicit expression for the isomorphism of
$(R\# A)\#\hat A$ to $R\ot (A\diamondsuit \hat A)$.
Of course, it would be possible to prove theorem 7.6 by using this map
and showing that it is actually an algebra isomorphism. However, it is
clear that the calculations will become very involved. On the other
hand, the technique we used here, by considering a faithful action of
our algebra, not only makes the proof easier, but it also shows where
the isomorphism comes from.
\snl
As a special case of this duality for actions of algebraic quantum
groups, we get of course the duality for actions of finite-dimensional
Hopf algebras ([vdB], see also [M2, page 167]).
Indeed, every finite-dimensional Hopf algebra is an
algebraic quantum group and the dual is precisely the dual Hopf algebra.
\nl
Now, we turn our attention again to the case of a general pairing. For
this, let us look again at the statement, preceding theorem 7.6. It says
that the algebra $P$ (which is isomorphic with the bismash product) is
spanned by the operators on $R\ot A$ of the form
$$x'\ot a' \to \sum ((S^{-1}a'_{(1)})x)x' \ot (a\# b) a'_{(2)}$$
where in the last factor, we use the natural left action of $B\# A$ on
$A$.
\snl
If we examine the above formula, we see that we have some kind of a
{\it twisting} of the algebras $R$ and $A\# B$. The action of $R$ on $R\ot A$
here is given by
$$x'\ot a' \to \sum ((S^{-1}a'_{(1)})x)x' \ot  a'_{(2)}.$$
In the case $(A,\hat A)$, it is essentially due to the fact that
$A\# \hat A$ acts by means of finite rank operators on $A$, that we are
able to {\it unfold} the effect of this twisting. In the general case we
realize this untwisting by a completely different technique. We add a
condition on the $A$-module $R$ to ensure that $Ax$ does not become too
big (finite-dimensional for the case where $A$ is an algebra with $1$)
for any $x\in R$.
\snl
We need the notion of a coaction (see [VD-Z3]).

\inspr{7.7} Definition \rm
If $R$ is an algebra over $\Bbb C$ and if $B$ is a regular multiplier
Hopf algebra, then a coaction of $B$ on $R$ is an injective homomorphism
$\Gamma : R \to M(R\ot B)$ such that $\Gamma(R)(1\ot B)\subseteq R\ot B$
and $(1\ot B)\Gamma(R)\subseteq R\ot B$
and such that $(\Gamma\ot\iota)\Gamma=(\iota\ot\Delta)\Gamma$.
\einspr

This last form of {\it coassociativity} can be given a meaning here in
different ways. We can e.g.\ write it as
$$(\Gamma\ot\iota)(\Gamma(x)(1\ot b))=(\iota\ot\Delta)(\Gamma(x))(1\ot
1\ot b).$$
Indeed, because $\Gamma(x)(1\ot b)\in R\ot B$, the left hand side makes
sense, while for the right hand side, just observe that $\Delta$ is
non-degenerate so that $\iota\ot \Delta$ can be extended to the
multiplier algebra $M(R\ot B)$.
\snl
It can be shown that for all coactions, we have the
equalities $\Gamma(R)(1\ot B)=R\ot B$
and $(1\ot B)\Gamma(R)=R\ot B$. In fact, the maps
$$\align x\ot b & \to \Gamma(x)(1\ot b) \\
         x\ot b & \to (1\ot b)\Gamma(x)
  \endalign$$
are injective on $R\ot B$ (see [VD-Z3]).
\snl
A trivial example of a coaction is the coaction of $B$ on itself,
given by the comultiplication $\Delta$. The bijectivity of the maps
above are in this case multiplier Hopf algebra axioms.
\snl
It is not hard to see that any coaction of  $B$ on $R$ will yield
an action of $A$ on $R$ by means of the formula
$$ax=(\iota\ot \omega_a)\Gamma(x)$$
where $\omega_a$ denotes the linear functional on $B$ given by
$b\to \langle a,b\rangle$. Remark, when $a=b'\trr a'$, we see that
$$(\iota\ot \omega_a)\Gamma(x)=(\iota\ot \omega_a')(\Gamma(x)(1\ot b'))$$
and this is well-defined. Hence, because the action of $B$ on $A$ is
unital, also $(\iota\ot \omega_a)\Gamma(x)$ makes sense.
Because of the fact that e.g.\ $\Gamma(R)(1\ot B)=R\ot B$, this
associated action is unital.
\snl
Not every action of $A$ will come from a coaction of $B$. However, if
$A$ is an algebraic quantum group and if $B=\hat A$, then this is the
case. So, for algebraic quantum groups, there is a one-to-one
correspondence between the actions of $A$ and the coactions of $\hat A$
(see [VD-Z3]).
\snl
Now, let us go back to the bismash product $(R\# A)\# B$ for a general
pair $(A, B)$. Assume that the action of $A$ on $R$ comes from a
coaction $\Gamma$ of $B$ on $R$. Let us use the Sweedler notation also
for these coactions. So, write $\Gamma(x)=\sum x_{(0)} \ot
x_{(1)}$ and make sure that $x_{(1)}$ is always covered by an element of
$B$. Then the algebra $P$ of operators on $R\ot A$
considered before, is spanned
by linear operators of the form
$$x'\ot a'\to \sum x_{(0)}x' \ot \langle S^{-1}a'_{(1)},x_{(1)}\rangle
(a\# b)a'_{(2)}.$$
In this formula, the element $a'_{(2)}$ is covered by $a\# b$ and so,
also $x_{(1)}$ will be covered through the pairing with $S^{-1}a'_{(1)}$
and this expression  has a meaning.
\snl
It is seen from the above formula that $P$ will be a subalgebra of the
tensor product algebra $R\ot Q$ where of course, $R$ acts on $R$ by left
multiplication and where $Q$ is the algebra generated by the left action of
$A\# B$ on $A$ on the one hand, and by the right action of $B$ on $A$ on
the other hand.
\snl
So, we clearly need some extra condition. It is similar to the
RL-condition as found in 9.4.5 of [M2]. In this context, we require that
the right action of an element $b\in B$ on $A$ is a multiplier of the
algebra $A\# B$  as acting on $A$. In other words, the algebra $Q$ above
coincides with the algebra $A\# B$ as acting on $A$.
\snl
There are some more arguments needed, but at least, this explains in
some sense the result that we get in [VD-Z3] and that is formulated in
the next theorem.

\inspr{7.8} Theorem \rm
If the action of $A$ on $R$ comes from a coaction of $B$ and if the
pairing between $A$ and $B$ is such that the right action of $B$ on $A$
is in the multiplier algebra of $A\# B$, as acting also on $A$, then,
$(R\# A)\# B$ is isomorphic with $R\ot (A\# B)$.
\einspr

This theorem is the natural generalization of the duality theorem for
actions of Hopf algebras as obtained by Blattner and Montgomery ([B-M],
[M2]) to the multiplier Hopf algebra case. The proof that we will give
in [VD-Z3] will use similar techniques as we have used to obtain theorem
7.6 and so, it is also slightly simpler that the one in [M2]. Also
remark that the extra condition on the pairing is always fulfilled for
the pair $(A,\hat A)$. Because also every action of $A$ comes from a
coaction of $\hat A$ in this case, theorem 7.8 also generalizes theorem
7.6.
\snl
The two duality theorems 7.6 and 7.8 can be combined with the fixed
point result in proposition 7.3 as follows. Consider the bidual action
of $A$ on the bismash product $(R\# A)\# B$. Use the isomorphism with
$R\ot (A\# B)$ to transform the bidual action to an action of $A$ on
$R\ot (A\# B)$. The multiplier algebra of $R\# A$ can be obtained as a
fixed point algebra. One could also try to use this idea to view the
smash product $R\# A$ itself as some kind of a fixed point algebra.
Similar techniques have been used in the theory of crossed products of
operator algebras by actions of locally compact groups to give such a
characterization of the crossed product in such a case (see e.g.\
theorem 3.11 in [VD1]).
%
%

\enddocument